\documentclass[12pt]{article}
\usepackage{amssymb}
\usepackage{amsmath}
\usepackage[doublespacing]{setspace}
\usepackage{geometry}
\usepackage[final,dvips]{graphicx}

\newtheorem{theorem}{Theorem}

\newtheorem{corollary}[theorem]{Corollary}

\newtheorem{example}{Example}

\newtheorem{lemma}[theorem]{Lemma}

\newtheorem{proposition}[theorem]{Proposition}

\input{tcilatex}
\begin{document}

\title{Signal Detection in High Dimension:\\
The Multispiked Case}
\author{ Alexei Onatski, Marcelo J. Moreira and Marc Hallin\thanks{%
Alexei Onatski: University of Cambridge, ao319@cam.ac.uk. Marcelo J.
Moreira: Funda\c{c}\~{a}o Getulio Vargas, mjmoreira@fgv.br. Marc Hallin:
Universit\'{e} libre de Bruxelles and Princeton University,
mhallin@ulb.ac.be.}}
\maketitle

\begin{abstract}
This paper deals with the local asymptotic structure, in the sense of Le
Cam's asymptotic theory of statistical experiments, of the signal detection
problem in high dimension. More precisely, we consider the problem of
testing the null hypothesis of sphericity of a high-dimensional covariance
matrix against an alternative of (unspecified) multiple symmetry-breaking
directions (\textit{multispiked} alternatives). Simple analytical
expressions for the asymptotic power envelope and the asymptotic powers of
previously proposed tests are derived. These asymptotic powers are shown to
lie very substantially below the envelope, at least for relatively small
values of the number of symmetry-breaking directions under the alternative.
In contrast, the asymptotic power of the likelihood ratio test based on the
eigenvalues of the sample covariance matrix is shown to be close to that
envelope. These results extend to the case of multispiked alternatives the
findings of an earlier study (Onatski, Moreira and Hallin, 2011) of the
single-spiked case. The methods we are using here, however, are entirely
new, as the Laplace approximations considered in the single-spiked context
do not extend to the multispiked case.
\end{abstract}

\bigskip

\section{Introduction}

In a recent paper, Onatski, Moreira and Hallin (2011) (hereafter OMH)
analyze the asymptotic power of statistical tests in the detection of a
signal in spherical real-valued Gaussian data as the dimensionality of the
data and the number of observations diverge to infinity at the same rate.
This paper generalizes OMH's alternative of a single symmetry-breaking
direction (\textit{single-spiked} alternative) to the alternative of
multiple symmetry-breaking directions (\textit{multispiked} alternative),
which is more relevant for applied work.

Contemporary tests of sphericity in a high-dimensional environment (see
Ledoit and Wolf (2002), Srivastava (2005), Schott (2006), Bai et al. (2009),
Chen et al. (2010), and Cai and Ma (2012)) consider general alternatives to
the null of sphericity. Our interest in alternatives with only a few
contaminating signals stems from the fact that in many applications, such as
speech recognition, macroeconomics, finance, wireless communication,
genetics, physics of mixture, and statistical learning, a few latent
variables typically explain a large portion of the variation in
high-dimensional data (see Baik and Silverstein (2006) for references). As a
possible explanation of this fact, Johnstone (2001) introduces the spiked
covariance model, where all eigenvalues of the population covariance matrix
of high-dimensional data are equal except for a small fixed number of
distinct \textquotedblleft spike eigenvalues.\textquotedblright\ The
alternative to the null of sphericity considered in this paper coincides
with Johnstone's model.

The extension from the single-spiked alternatives\ of OMH to the
multi-spiked alternatives\ considered here, however, is all but
straightforward. The difficulty arises because the extension of the main
technical tool in OMH (Lemma~2), which analyzes high-dimensional spherical
integrals, to integrals over high-dimensional real Stiefel manifolds
obtained in Onatski (2012) is not easily amenable to the Laplace
approximation method used in OMH. Therefore, in this paper, we develop a
completely different technique, inspired from the large deviation analysis
of spherical integrals by Guionnet and Maida (2005).

Let us describe the setting and main results in more detail. Suppose that
the data consist of $n$ independent observations $X_{t},$ $t=1,...,n$ of a $%
p $-dimensional Gaussian vector with mean zero and positive definite
covariance matrix~$\Sigma $. Let $\Sigma =\sigma ^{2}\left(
I_{p}+VHV^{\prime }\right) ,$ where $I_{p}$ is the $p$-dimensional identity
matrix, $\sigma $ is a scalar, $H$ an $r\times r$ diagonal matrix with
elements $h_{j}\geq 0,$ $j=1,...,r$ along the diagonal, and $V$ a $\left(
p\times r\right) $-dimensional parameter normalized so that $V^{\prime
}V=I_{r}$. We are interested in the asymptotic power of tests of the null
hypothesis $H_{0}:h_{1}=...=h_{r}=0$ against the alternative $H_{1}:$ $%
h_{j}>0$ for some $j=1,...,r,$ based on the eigenvalues of the sample
covariance matrix of the data when $n,p\rightarrow \infty $ so that $%
p/n\rightarrow c$ with $0<c<\infty $, an asymptotic regime which we
abbreviate into $n,p\rightarrow _{c}\infty $. The matrix $V$ is an
unspecified nuisance parameter, the columns of which indicate the directions
of the perturbations of sphericity.

We consider the cases of specified and unspecified $\sigma ^{2}.$ For the
sake of simplicity, in the rest of this introduction, we only discuss the
case of specified $\sigma ^{2}=1,$ although the case of unspecified $\sigma
^{2}$ is more realistic. Denoting by $\lambda _{j}$ the $j$-th largest
sample covariance eigenvalue, let $\lambda =\left( \lambda _{1},...,\lambda
_{m}\right) ,$ where $m=\min \left( n,p\right) $. We begin our analysis with
a study of the asymptotic properties of the likelihood ratio process $%
\left\{ L\left( h;\lambda \right); h\in \left[ 0,\bar{h}\right] ^{r}\right\}
,$ where $h=\left( h_{1},...,h_{r}\right) ,$ $\bar{h}\in \left[ 0,\sqrt{c}%
\right) $ and $L\left( h;\lambda \right) $ is defined as the ratio of the
density of $\lambda $ under $H_{1}$ to that under $H_{0}$, considered as a $%
\lambda $-measurable random variable. Note that $L\left( h;\lambda \right) $
depends on $n$ and $p,$ while $\lambda $ is $m=\min \left\{ n,p\right\} $
-dimensional. An exact formula for $L\left( h;\lambda \right) $ involves the
integral $\int_{\mathcal{O}\left( p\right) }e^{\limfunc{tr}\left(
AQBQ^{\prime }\right) }\left( dQ\right) $ over the orthogonal group $%
\mathcal{O}\left( p\right) $, where the $p\times p$ matrix $A$ has a
deficient rank $r$. In the single-spiked case ($r=1$), OMH link this
integral to the confluent form of the Lauricella function, and use this link
to establish a representation of the integral in the form of a contour
integral (see Wang (2010) and Mo (2011) for independent different
derivations of this contour integral representation for this particular $r=1$
case). Then, the Laplace approximation to the contour integral is used to
derive the asymptotic behavior of $L\left( h;\lambda \right) $.

Onatski (2012) generalizes the contour integral representation to the
multispiked case ($r>1$). For complex-valued data, such a generalization
allows him to extend OMH's results to the multi-spiked context.
Unfortunately, for real-valued data, which we are concerned with in this
paper, this generalization is not straightforwardly amenable to the Laplace
approximation method. Therefore, in this paper, we consider a totally
different approach. For the $r=1$ case, Guionnet and Maida (2005) (hereafter
GM) use large deviation methods to derive a second-order asymptotic
expansion of $\int_{\mathcal{O}\left( p\right) }e^{\limfunc{tr}\left(
AQBQ^{\prime }\right) }\left( dQ\right) $ as the non-zero eigenvalues of $A$
diverge to infinity (see their Theorem 3). We extend GM's second-order
expansion to the $r>1$ case, and use that extension to derive the
asymptotics of $L\left( h;\lambda \right) $.

More precisely, we show that, for any $\bar{h}$ such that $0<\bar{h}<\sqrt{c}%
,$ the sequence of log-likelihood processes $\{\ln L\left( h;\lambda \right)
;h\in \lbrack 0,\bar{h}]^{r}\}$ converges weakly to a Gaussian process $%
\left\{ \mathcal{L}_{\lambda }(h);h\in \left[ 0,\bar{h}\right] ^{r}\right\} $
under the null hypothesis $H_{0}$ as $n,p\rightarrow _{c}\infty $. The index 
$\lambda $ in the notation $\mathcal{L}_{\lambda }(h)$ is used to
distinguish the limiting $\lambda $-log-likelihood process in the case of
specified $\sigma ^{2}=1,$ from that of the $\mu $-log-likelihood process
considered in the case of unspecified $\sigma ^{2},$ which we denote by $%
\mathcal{L}_{\mu }(h)$ (see Section~2). The limiting process has mean $%
\limfunc{E}\left[ \mathcal{L}_{\lambda }\mathcal{(}h\mathcal{)}\right] =%
\frac{1}{4}\sum\limits_{i,j=1}^{r}\ln \left( 1-h_{i}h_{j}/c\right) $ and
autocovariance function $\limfunc{Cov}\left( \mathcal{L}_{\lambda }\left(
h\right) ,\mathcal{L}_{\lambda }(\tilde{h})\right) =-\frac{1}{2}%
\sum_{i,j=1}^{r}\ln \left( 1-h_{i}\tilde{h}_{j}/c\right) .$ That convergence
entails the weak convergence, in the Le Cam sense, of the $h$-indexed
statistical experiments $\mathcal{E}_{\lambda }^{m}$ under which the
eigenvalues $\lambda _{1},...,\lambda _{m}$ are observed, i.e. the
statistical experiments with log-likelihood process $\left\{ \ln L_{\lambda
}\left( h\right) ;h\in \left[ 0,\bar{h}\right] ^{r}\right\} $ (see van der
Vaart (1998), page 126). Although this limiting process is Gaussian, it is
not a log-likelihood process of the Gaussian shift type, so that the
statistical experiments $\mathcal{E}_{\lambda }^{m}$ under study are not
locally asymptotically normal (LAN) ones. The weak convergence of $\mathcal{E%
}_{\lambda }^{m}$ implies, however, via Le Cam's first lemma (see van der
Vaart 1998, p.88), that the joint distributions of the normalized sample
covariance eigenvalues under the null and under alternatives associated with 
$h\in \left[ 0,\sqrt{c}\right) $ are mutually contiguous.

An asymptotic power envelope for $\lambda $-based tests of $H_{0}$ against $%
H_{1}$ can be constructed using the Neyman-Pearson lemma and Le Cam's third
lemma. We show that, for tests of size $\alpha ,$ the maximum achievable
asymptotic power against a point alternative $h=(h_{1},...,h_{r})$ equals $%
1-\Phi \left[ \Phi ^{-1}\left( 1-\alpha \right) -\sqrt{W(h)}\right] ,$ where 
$\Phi $ is the standard normal distribution function and $W(h)=-\frac{1}{2}%
\sum_{i,j=1}^{r}\ln \left( 1-h_{i}h_{j}/c\right) $. As we explain in the
paper, this asymptotic power envelope is valid not only for the $\lambda $%
-based tests, but also for all tests that are invariant under left orthogonal
transformations of the data $X_{t},$ $t=1,...,n$.

Next, we consider previously proposed tests of sphericity and of the
equality of the population covariance matrix to a given matrix . We focus on
the tests studied in Ledoit and Wolf (2002), Bai et al (2009), and Cai and
Ma (2012). We find that, in general, the asymptotic powers of those tests
are substantially lower than the corresponding asymptotic power envelope
value. In contrast, our computations for the case $r=2$ show that the
asymptotic powers of the $\lambda $- and $\mu $-based likelihood ratio tests
are close to the power envelope.

The rest of the paper is organized as follows. Section 2 establishes the
weak convergence of the log-likelihood ratio process to a Gaussian process.
Section 3 provides an analysis of the asymptotic powers of various
sphericity tests, derives the asymptotic power envelope, and proves its
validity for general invariant tests. Section 4 concludes. All proofs are
given in the Appendix.

\section{Asymptotics of likelihood ratio processes}

Let $X$ be a $p\times n_{p}$ matrix with independent Gaussian $N\left(
0,\sigma ^{2}\left( I_{p}+VHV^{\prime }\right) \right) $ columns. Let $%
\lambda _{p1}\geq ...\geq \lambda _{pp}$ be the ordered eigenvalues of $%
\frac{1}{n_{p}}XX^{\prime }$ and write $\lambda _{p}=\left( \lambda
_{p1},...,\lambda _{pm}\right) ,$ where $m\!=\!\min \left\{
p,\!n_{p}\right\} $. Similarly, let $\mu _{pi}=\!\lambda _{pi}/\!\left(
\lambda _{p1}\!+\!...\!+\!\lambda _{pp}\right) ,$ $i=1,...,m$ and $\mu
_{p}\!=\!\left( \mu _{p1},...,\mu _{p,m-1}\right) $.

As explained in the introduction, our goal is to study the asymptotic power,
as $n_{p},p\rightarrow _{c}\infty ,$ of the eigenvalue-based tests of $%
H_{0}:h_{1}=...=h_{r}=0$ against $H_{1}:h_{j}>0$ for some $i=1,...,r$, where 
$h_{j}$ are the diagonal elements of the diagonal matrix $H$. If $\sigma ^{2}
$ is specified, the model is invariant with respect to left and right
orthogonal transformations; sufficiency and invariance arguments (see
Appendix 5.4 for details) lead to considering tests based on $\lambda _{p}$
only. If~$\sigma ^{2}$ is unspecified, the model is invariant with respect
to left and right orthogonal transformations and multiplications by non-zero
scalars; sufficiency and invariance arguments (see Appendix 5.4) lead to
considering tests based on $\mu _{p}$ only. Note that the distribution of $%
\mu _{p}$ does not depend on $\sigma ^{2},$ whereas, if $\sigma ^{2}$ is
specified, we can always normalize $\lambda _{p}$ dividing it by $\sigma
^{2}.$ Therefore, we henceforth assume without loss of generality that~$%
\sigma ^{2}=1$.

Let us denote the joint density of $\lambda _{p1},...,\lambda _{pm}$ at $%
\tilde{x}=\left( x_{1},...,x_{m}\right) \in (\mathbb{R}^{+})^{m}$ as~$%
f_{\lambda p}\left( \tilde{x};h\right) $, and that of $\mu _{p1},...,\mu
_{p,m-1}$ at $\tilde{y}=\left( y_{1},...,y_{m-1}\right) \in (\mathbb{R}%
^{+})^{m-1}$ as $f_{\mu p}\left( \tilde{y};h\right) $. We have%
\begin{equation}
f_{\lambda p}\left( \tilde{x};h\right) =\tilde{\gamma}\frac{%
\prod_{i=1}^{m}x_{i}^{\frac{\left\vert p-n_{p}\right\vert -1}{2}%
}\prod_{i<j}^{m}\left( x_{i}-x_{j}\right) }{\prod_{j=1}^{r}\left(
1+h_{j}\right) ^{n_{p}/2}}\int\limits_{\mathcal{O}\left( p\right) }e^{-\frac{%
n_{p}}{2}\limfunc{tr}\left( \Pi Q^{\prime }\mathcal{X}Q\right) }\left( 
\mathrm{d}Q\right) ,  \label{common complex1}
\end{equation}%
where $\tilde{\gamma}$ depends only on $n_{p}$ and $p$; $\Pi =\mathrm{diag}%
\left( \left( 1+h_{1}\right) ^{-1},...,\left( 1+h_{r}\right)
^{-1},1,...,1\right) $; $\mathcal{X}=\mathrm{diag}\left(
x_{1},...,x_{m},0,...,0\right) $ is a $\left( p\times p\right) $ diagonal
matrix; $\mathcal{O}\left( p\right) $ is the set of all $p\times p$
orthogonal matrices; and $\left( \mathrm{d}Q\right) $ is the invariant
measure on the orthogonal group~$\mathcal{O}\left( p\right) $, normalized to
make the total measure unity. Formula (\ref{common complex1}) is a special
case of the density given in James (1964, p.483) for $n_{p}\geq p,$ and
follows from Theorems 2 and 6 in Uhlig (1994) for $n_{p}<p$.

Let $x=x_{1}+...+x_{m}$ and let $y_{i}=x_{i}/x.$ Note that the Jacobian of
the coordinate change from $\left( x_{1},...,x_{m}\right) $ to $\left(
y_{1},...,y_{m-1},x\right) $ is $x^{m-1}.$ Changing variables in (\ref%
{common complex1}) and integrating $x$ out, we obtain%
\begin{equation}
f_{\mu p}\left( \tilde{y};h\right) =\tilde{\gamma}\frac{%
\prod_{i=1}^{m}y_{i}^{\frac{\left\vert p-n_{p}\right\vert -1}{2}%
}\prod_{i<j}^{m}\left( y_{i}-y_{j}\right) }{\prod_{j=1}^{r}\left(
1+h_{j}\right) ^{n_{p}/2}}\int_{0}^{\infty }x^{\frac{n_{p}p}{2}%
-1}\int\limits_{\mathcal{O}\left( p\right) }e^{-\frac{n_{p}}{2}x\limfunc{tr}%
\left( \Pi Q^{\prime }\mathcal{Y}Q\right) }\left( \mathrm{d}Q\right) \mathrm{%
d}x,  \label{common complex2}
\end{equation}%
where $\mathcal{Y}=\mathrm{diag}\left( y_{1},...,y_{m},0,...,0\right) $ is a 
$\left( p\times p\right) $ diagonal matrix.

Consider the likelihood ratios $L_{p}\left( h;\lambda _{p}\right)
=f_{\lambda p}\left( \lambda _{p};h\right) /f_{\lambda p}\left( \lambda
_{p};0\right) $ and $L_{p}\left( h;\mu \right) =f_{\mu p}\left( \mu
_{p};h\right) /f_{\mu p}\left( \mu _{p};0\right) $. Formulae (\ref{common
complex1}) and (\ref{common complex2}) imply the following proposition.

\begin{proposition}
\label{Proposition1}\textit{Let }$\mathcal{O}\left( p\right) $ \textit{be
the set of all }$p\times p$\textit{\ orthogonal matrices. Denote by }$\left( 
\mathrm{d}Q\right) $\textit{\ the invariant measure on the orthogonal group} 
$\mathcal{O}\left( p\right) $ \textit{normalized to make the total measure
unity. Put }$\Lambda _{p}=\mathrm{diag}\left( \lambda _{p1},...,\lambda
_{pp}\right) ,$ $S_{p}=\lambda _{p1}+...+\lambda _{pp},$ and let $D_{p}$ be
the $p\times p$ diagonal matrix $\mathrm{diag}\left( \frac{1}{2c_{p}}\frac{%
h_{1}}{1+h_{1}},...,\frac{1}{2c_{p}}\frac{h_{r}}{1+h_{r}},0,...,0\right) ,$
where $c_{p}=p/n_{p}$. \textit{Then,}%
\begin{eqnarray}
L_{p}\!\left( h;\lambda _{p}\right) \!\!\!\!
&=&\!\!\!\!\prod_{j=1}^{r}\left( 1+h_{j}\right) ^{-\frac{n_{p}}{2}%
}\int\limits_{\mathcal{O}\left( p\right) }e^{p\limfunc{tr}\left(
D_{p}Q^{\prime }\Lambda _{p}Q\right) }\left( \mathrm{d}Q\right) \text{ 
\textit{and}}  \label{LR1} \\
L_{p}\!\left( h;\mu _{p}\right) \!\!\!\! &=&\!\!\!\!\prod_{j=1}^{r}\left(
1+h_{j}\right) ^{-\frac{n_{p}}{2}}\!\frac{\left( \frac{n_{p}}{2}\right) ^{%
\frac{n_{p}p}{2}}}{\Gamma \!\left( \frac{n_{p}p}{2}\right) }%
\!\int_{0}^{\infty }\!\!\!x^{\frac{n_{p}p}{2}-1}e^{-\frac{n_{p}}{2}%
x}\!\!\int\limits_{\mathcal{O}\left( p\right) }\!\!\!e^{p\frac{x}{S_{p}}%
\limfunc{tr}\left( D_{p}Q^{\prime }\Lambda _{p}Q\right) }\left( \!\mathrm{d}%
Q\!\right) \!\mathrm{d}x.  \label{LR2}
\end{eqnarray}
\end{proposition}

In the special case where $r=1,$ the rank of the matrix $D_{p}$ equals one,
and the integrals over the orthogonal group in (\ref{LR1}) and (\ref{LR2})
can be rewritten as integrals over a $p$-dimensional sphere. OMH show how
such spherical integrals can be represented in the form of contour
integrals, and apply Laplace approximation to these contour integrals to
establish the asymptotic properties of $L_{p}\left( h;\lambda _{p}\right) $
and $L_{p}\left( h;\mu _{p}\right) .$ In the~$r>1$ case, the integrals in (%
\ref{LR1}) and (\ref{LR2}) can be rewritten as integrals over a Stiefel
manifold, the set of all orthonormal $r$\textit{-frames} in $\mathbb{R}^{p}$%
. Onatski (2012) obtains a generalization of the contour integral
representation from spherical integrals to integrals over Stiefel manifolds.
Unfortunately, the Laplace approximation method does not straightforwardly
extend to that generalization, and we therefore propose an alternative
method of analysis.

The second-order asymptotic behavior, as $p$ goes to infinity, of integrals
of the form $\int_{\mathcal{O}\left( p\right) }e^{p\limfunc{tr}\left(
DQ^{\prime }\Lambda Q\right) }\left( \mathrm{d}Q\right) $ was analyzed in
Guionnet and Maida (2005) (Theorem~3) for the particular case where $D$ is a
fixed matrix of rank one, $\Lambda $ a deterministic matrix, and under the
condition that the empirical distribution of $\Lambda $'s eigenvalues
converges to a distribution function with bounded support. Below, we extend
Guionnet and Maida's approach to cases where $D=D_{p}$ has rank larger than
one, and to the stochastic setting of this paper. We then use such an
extension to derive the asymptotic properties of $L_{p}\left( h;\lambda
_{p}\right) $ and $L_{p}\left( h;\mu _{p}\right) $.

Let $\hat{F}_{p}^{\lambda }$ be the empirical distribution of $\lambda
_{p1},...,\lambda _{pp}$, and denote by $F_{p}^{MP}$ the Marchenko-Pastur
distribution function, with density%
\begin{equation}
f_{p}^{MP}\left( x\right) =\frac{1}{2\pi c_{p}x}\sqrt{\left( b_{p}-x\right)
\left( x-a_{p}\right) },  \label{MP density}
\end{equation}%
where $a_{p}=\left( 1-\sqrt{c_{p}}\right) ^{2}$ and $b_{p}=\left( 1+\sqrt{%
c_{p}}\right) ^{2}$, and a mass of $\max \left( 0,1-c_{p}^{-1}\right) $ at
zero. As is well known, the difference between $\hat{F}_{p}^{\lambda }$ and $%
F_{p}^{MP}$ weakly converges to zero a.s. as $p,n_{p}\rightarrow _{c}\infty $%
. Moreover, $\lambda _{p1}\overset{a.s}{\rightarrow }\left( 1+\sqrt{c}%
\right) ^{2},$ and $\lambda _{pp}\overset{a.s}{\rightarrow }\left( 1-\sqrt{c}%
\right) ^{2}$ if $c>1$, and $\lambda _{pp}\overset{a.s}{\rightarrow }0$ if $%
c\leq 1.$

Consider the Hilbert transform of $F_{p}^{MP},$ $H_{p}^{MP}(x)=\int \left(
x-\lambda \right) ^{-1}\mathrm{d}F_{p}^{MP}\left( \lambda \right) .$ That
transform is well defined for real $x$ outside the support of $F_{p}^{MP}$,
that is, on the set $\mathbb{R}\backslash \limfunc{supp}\left(
F_{p}^{MP}\right) .$ Using (\ref{MP density}), we get%
\begin{equation}
H_{p}^{MP}\left( x\right) =\frac{x+c_{p}-1-\sqrt{\left( x-c_{p}-1\right)
^{2}-4c_{p}}}{2c_{p}x},  \label{Stijeltjes analytic}
\end{equation}%
where the sign of the square root is chosen to be the sign of $\left(
x-c_{p}-1\right) $. It is not hard to see that $H_{p}^{MP}\left( x\right) $
is strictly decreasing on $\mathbb{R}\backslash \limfunc{supp}\left(
F_{p}^{MP}\right) $. Thus, on $H_{p}^{MP}\left( \mathbb{R}\backslash 
\limfunc{supp}\left( F_{p}^{MP}\right) \right) $, we can define an inverse
function $K_{p}^{MP}$, with values%
\begin{equation}
K_{p}^{MP}\left( x\right) =\frac{1}{x}+\frac{1}{1-c_{p}x},\text{ }x\in
H_{p}^{MP}\left( \mathbb{R}\backslash \limfunc{supp}\left( F_{p}^{MP}\right)
\right) .  \label{K transform}
\end{equation}%
The so-called $R$-transform $R_{p}^{MP}$ of $F_{p}^{MP}$ takes the form 
\begin{equation*}
R_{p}^{MP}\left( x\right) =K_{p}^{MP}\left( x\right) -1/x=1/\left(
1-c_{p}x\right) .
\end{equation*}

For $\varepsilon >0$ and $\eta >0$ sufficiently small, consider the subset
of $\mathbb{R}$ 
\begin{equation*}
\Omega _{\varepsilon \eta }=\left\{ \!\!\!%
\begin{array}{lc}
\left[ -\eta ^{-1},0\right) \cup \left( 0,\frac{1}{\sqrt{c}\left( 1+\sqrt{c}%
\right) }-\varepsilon \right] & \text{for }c\geq 1, \\ 
\left[ -\frac{1}{\sqrt{c}\left( 1-\sqrt{c}\right) }+\varepsilon ,0\right)
\cup \left( 0,\frac{1}{\sqrt{c}\left( 1+\sqrt{c}\right) }-\varepsilon \right]
& \text{for }c<1.%
\end{array}%
\right.
\end{equation*}%
From (\ref{Stijeltjes analytic}), $H_{p}^{MP}\left( \mathbb{R}\backslash 
\limfunc{supp}\left( F_{p}^{MP}\right) \right) =\left( -\infty ,0\right)
\cup \left( 0,\frac{1}{\sqrt{c_{p}}\left( 1+\sqrt{c_{p}}\right) }\right)
\cup \left( \frac{1}{\sqrt{c_{p}}\left( \sqrt{c_{p}}-1\right) },\infty
\right) $ when $c_{p}>1$, $\left( -\frac{1}{\sqrt{c_{p}}\left( 1-\sqrt{c_{p}}%
\right) },0\right) \cup \left( 0,\frac{1}{\sqrt{c_{p}}\left( 1+\sqrt{c_{p}}%
\right) }\right) $ when $c_{p}<1$, and $\left( -\infty ,0\right) \cup \left(
0,1/2\right) $ when $c_{p}=1$. Therefore, $\Omega _{\varepsilon \eta
}\subset H_{p}^{MP}\left( \mathbb{R}\backslash \limfunc{supp}\left(
F_{p}^{MP}\right) \right) $ with probability approaching one as $%
n_{p},p\rightarrow _{c}\infty .$

\begin{proposition}
\label{Proposition2}Let $\left\{ \Theta _{p}\right\} $ be a sequence of
random $p\times p$ diagonal matrices \linebreak $\limfunc{diag}\left( \theta
_{p1},...,\theta _{pr},0,...,0\right) ,$ where $\theta _{pj}\neq 0,$ $%
j=1,...,r.$ Further, let $v_{pj}=R_{p}^{MP}\left( 2\theta _{pj}\right) $,
where $R_{p}^{MP}\left( x\right) =1/\left( 1-c_{p}x\right) $ is the $R$%
-transform of the Marchenko-Pastur distribution $F_{p}^{MP}$. Assume that,
for some $\varepsilon >0$ and $\eta >0,$ $2\theta _{pj}\in \Omega
_{\varepsilon ,\eta }$ with probability approaching one as $%
n_{p},p\rightarrow _{c}\infty $. Then,%
\begin{eqnarray*}
\int\limits_{\mathcal{O}\left( p\right) }e^{p\limfunc{tr}\left( \Theta
_{p}Q^{\prime }\Lambda _{p}Q\right) }\left( \mathrm{d}Q\right)
&=&e^{p\sum_{j=1}^{r}\left[ \theta _{pj}v_{pj}-\frac{1}{2p}\sum_{i=1}^{p}\ln
\left( 1+2\theta _{pj}v_{pj}-2\theta _{pj}\lambda _{p,i}\right) \right] } \\
&&\times \prod_{j=1}^{r}\prod_{s=1}^{j}\sqrt{1-4\left( \theta
_{pj}v_{pj}\right) \left( \theta _{ps}v_{ps}\right) c_{p}}\left(
1+o(1)\right) \text{ a.s.},
\end{eqnarray*}%
where $o(1)$ is uniform over all sequences $\left\{ \Theta _{p}\right\} $
satisfying the assumption.
\end{proposition}

This proposition extends Theorem 3 of Guionnet and Maida (2005) to cases
when $\limfunc{rank}\left( \Theta _{p}\right) >1,$ $\theta _{pj}$ depends on 
$p,$ and $\Lambda _{p}$ is random. When $r=1,$ $\theta _{p1}=\theta >0$ and $%
v_{p1}=v$ are fixed, it is straightforward to verify that $\sqrt{1-4\theta
^{2}v^{2}c_{p}}=\sqrt{4\theta ^{2}}/\sqrt{Z},$ where $Z=\int \left(
K_{p}^{MP}\left( 2\theta \right) -\lambda \right) ^{-2}\mathrm{d}%
F_{p}^{MP}\left( \lambda \right) .$ In Guionnet and Maida's (2005) Theorem
3, the expression $\sqrt{4\theta ^{2}}/\sqrt{Z}$ should have been used
instead of $\sqrt{Z-4\theta ^{2}}/\theta \sqrt{Z},$ which is a typo.

Setting $r=1$ and $\theta _{p1}=\frac{1}{2c_{p}}\frac{h}{1+h}$ in
Proposition \ref{Proposition2} and using formula~(\ref{LR1}) from
Proposition \ref{Proposition1} gives us an expression for $L_{p}\!\left(
h;\lambda _{p}\right) $ which is an equivalent of formula (4.1) in Theorem 7
of OMH. Theorem \ref{Theorem1} below uses Proposition \ref{Proposition2} to
generalize Theorem 7 of OMH to the multispiked case $r>1$.

Let $\theta _{pj}=h_{j}/2c_{p}\left( 1+h_{j}\right) $ and%
\begin{equation}
H_{\delta }\!=\!\left\{ \!\!\!%
\begin{array}{lc}
\left[ -1+\delta ,0\right) \cup \left( 0,\sqrt{c}-\delta \right] & \text{for 
}c>1, \\ 
\left[ -\sqrt{c}+\delta ,0\right) \cup \left( 0,\sqrt{c}-\delta \right] & 
\text{for }c\leq 1.%
\end{array}%
\right.  \label{h domain}
\end{equation}%
The condition $h_{j}\in H_{\delta }$ for some $\delta >0$ implies that $%
2\theta _{pj}\in \Theta _{\varepsilon \eta }$ for some $\varepsilon >0,$ $%
\eta >0$ and $p$ sufficiently large. Below, we are only interested in
non-negative values of $h_{j}$, and assume that $h_{j}\in \left( 0,\sqrt{c}%
-\delta \right] $ under the alternative hypothesis. The corresponding $%
\theta _{pj}$, thus, is positive.

With the above setting for $\theta _{pj},$ we have $v_{pj}=1+h_{j}$ and $%
K_{p}^{MP}\left( 2\theta _{pj}\right) =\left( c_{p}+h_{j}\right) \left(
1+h_{j}\right) /h_{j}=z_{j0},$ say, as in Theorem 7 in OMH. Define\vspace{%
-0.2in}%
\begin{equation}
\Delta _{p}\left( z_{j0}\right) =\sum_{i=1}^{p}\ln \left( z_{j0}-\lambda
_{pi}\right) -p\int \ln \left( z_{j0}-\lambda \right) \mathrm{d}%
F_{p}^{MP}\left( \lambda \right) .  \label{delta_definition}
\end{equation}

\begin{theorem}
\label{Theorem1}\textit{Suppose that the null hypothesis is true (}$h=0$)%
\textit{. Let }$\delta $\textit{\ be any fixed number such that }$0<\delta <%
\sqrt{c}$, \textit{and let }$C\left[ 0,\sqrt{c}-\delta \right] ^{r}$\textit{%
\ be the space of real-valued continuous functions on }$\left[ 0,\sqrt{c}%
-\delta \right] ^{r}$\textit{\ equipped with the supremum norm}. \textit{%
Then,} \textit{as }$p,n_{p}\rightarrow _{c}\infty ,$\textit{\ }%
\begin{eqnarray}
L_{p}\!\left( h;\lambda _{p}\right) \!\!\!\!
&=&\!\!\!\!\prod\limits_{j=1}^{r}\!\exp \!\left\{ \!-\frac{1}{2}\Delta
_{p}\left( z_{j0}\right) \!+\!\frac{1}{2}\sum\limits_{s=1}^{j}\ln \!\left(
\!1\!-\!\frac{h_{j}h_{s}}{c_{p}}\!\right) \!\right\} \left( 1\!+\!o\left(
1\right) \right) \text{ \textit{and}}  \label{equivalence 1} \\
L_{p}\!\left( h;\mu _{p}\!\right) \!\!\!\! &=&\!\!\!\!L_{p}\!\left(
h;\lambda _{p}\right) \exp \!\left\{ \!\frac{1}{4c_{p}}\left(
\!\sum_{j=1}^{r}\!h_{j}\!\right) ^{2}\!-\frac{S_{p}\!-\!p}{2c_{p}}%
\sum_{j=1}^{r}\!h_{j}\!\right\} \!\left( 1\!+\!o\left( 1\right) \right) \!,
\label{equivalence 2}
\end{eqnarray}%
\textit{almost surely, }where the $o\left( 1\right) $ terms are uniform in $%
h\in \left[ 0,\sqrt{c}-\delta \right] ^{r}$. \textit{Furthermore, }$\ln
L_{p}\!\left( h;\lambda _{p}\right) $ \textit{and} $\ln L_{p}\!\left( h;\mu
_{p}\!\right) ,$\textit{\ viewed as random elements of }$C\left[ 0,\sqrt{c}%
-\delta \right] ^{r}$,\textit{\ converge\ weakly to }$\mathcal{L}_{\lambda
}\left( h\right) $ \textit{and} $\mathcal{L}_{\mu }\left( h\right) $ \textit{%
with Gaussian finite-dimensional distributions such that }$\limfunc{E}\left( 
\mathcal{L}_{\lambda }\left( h\right) \right) =-\frac{1}{2}\limfunc{Var}%
\left( \mathcal{L}_{\lambda }\left( h\right) \right) $\textit{, }$\limfunc{E}%
\left( \mathcal{L}_{\mu }\left( h\right) \right) =-\frac{1}{2}\limfunc{Var}%
\left( \mathcal{L}_{\mu }\left( h\right) \right) ,$ and, \textit{for any }$h,%
\tilde{h}\in \left[ 0,\sqrt{c}-\delta \right] ^{r},$%
\begin{eqnarray}
&&\limfunc{Cov}\left( \mathcal{L}_{\lambda }\left( h\right) ,\mathcal{L}%
_{\lambda }\left( \tilde{h}\right) \right) =-\frac{1}{2}\sum_{i,j=1}^{r}\ln
\left( 1-\frac{h_{i}\tilde{h}_{j}}{c}\right) \text{, and}  \label{CovLambda}
\\
&&\limfunc{Cov}\left( \mathcal{L}_{\mu }\left( h\right) ,\mathcal{L}_{\mu
}\left( \tilde{h}\right) \right) =-\frac{1}{2}\sum_{i,j=1}^{r}\left( \ln
\left( 1-\frac{h_{i}\tilde{h}_{j}}{c}\right) +\frac{h_{i}\tilde{h}_{j}}{c}%
\right) \text{.}  \label{CovMu}
\end{eqnarray}
\end{theorem}

Theorem \ref{Theorem1} and Le Cam's first lemma (van der Vaart (1998), p.88)
imply that the joint distributions of $\lambda _{1},...,\lambda _{m}$ (as
well as those of $\mu _{1},...,\mu _{m-1}$) under the null and under the
alternative are mutually contiguous for any $h\in \left[ 0,\sqrt{c}\right)
^{r}$. By applying Le Cam's third lemma (van der Vaart (1998), p.90), we can
study the \textquotedblleft local\textquotedblright\ powers of tests
detecting signals in noise. The requirement that $h_{j}$ be positive under
alternatives corresponds to situations where the signals contained in the
data are independent from the noise. If dependence between the signals and
the noise is allowed, one might consider two-sided alternatives of the form $%
H_{1}:h_{j}\neq 0$ for some $j$. Values of $h_{j}$ between $-1$ and $0$
correspond to alternatives under which the noise variance is reduced along
certain directions. In view of Proposition \ref{Proposition2}, it should not
be difficult to generalize Theorem \ref{Theorem1} to the case of fully ($%
h_{j}\neq 0,$ all $j$'s) or partially ($h_{j}\neq 0,$ some $j$'s) two-sided
alternatives. This problem will not be discussed here, and is left for
future research.

\section{Asymptotic power analysis}

Denote by $\beta _{\lambda }\left( h\right) $ and $\beta _{\mu }\left(
h\right) ,$ respectively, the asymptotic powers of the asymptotically most
powerful $\lambda $- and $\mu $-based tests of size $\alpha $ of the null $%
h=0$ against a point alternative $h=\left( h_{1},...,h_{r}\right) \neq 0$
with $h_{j}<\sqrt{c},$ $j=1,...,r$. As functions of $h,$ $\beta _{\lambda }$
and $\beta _{\mu }$ are called the \textit{asymptotic power envelopes}.

\begin{proposition}
\label{Proposition3}\textit{Let }$\Phi $\textit{\ denote the standard normal
distribution function. Then,}%
\begin{eqnarray}
\beta _{\lambda }\left( h\right) \! &=&\!1\!-\!\Phi \left[ \!\Phi
^{-1}\left( 1\!-\!\alpha \right) \!-\!\sqrt{-\frac{1}{2}\sum_{i,j=1}^{r}\ln
\left( 1\!-\!\frac{h_{i}h_{j}}{c}\right) }\,\right] \text{ \textit{and}}
\label{local power} \\
\beta _{\mu }\left( h\right) \! &=&\!1\!-\!\Phi \left[ \!\Phi ^{-1}\left(
1\!-\!\alpha \right) \!-\!\sqrt{-\frac{1}{2}\sum_{i,j=1}^{r}\left( \ln
\left( 1\!-\!\frac{h_{i}h_{j}}{c}\right) \!+\!\frac{h_{i}h_{j}}{c}\right) }\,%
\right] .  \label{local power mu}
\end{eqnarray}
\end{proposition}

\begin{figure}[t!]
\centering
\includegraphics[width=4in]{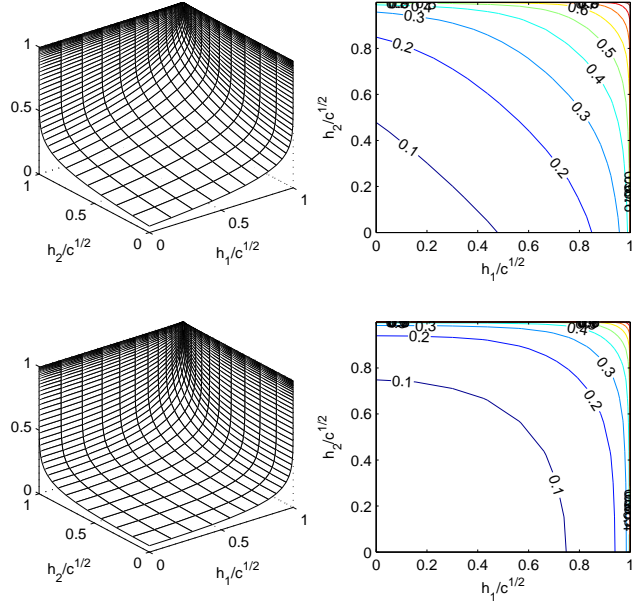}
\caption{The power envelopes $\protect%
\beta _{\protect\lambda }\left( h\right) $ (upper panel) and $\protect\beta %
_{\protect\mu }\left( h\right) $ (lower panel) for $\protect\alpha =0.05,$
as functions of $h/\protect\sqrt{c}=\left( h_{1},h_{2}\right) /\protect\sqrt{%
c}$.}
\label{envelope_may2012}
\end{figure}

Figure \ref{envelope_may2012} shows the asymptotic power envelopes $\beta
_{\lambda }\left( h\right) $ and $\beta _{\mu }\left( h\right) $ as
functions of $h_{1}/\sqrt{c}$ and $h_{2}/\sqrt{c}$ when $h=\left(
h_{1},h_{2}\right) $ is two-dimensional.


It is important to realize that the asymptotic power envelopes derived in
Proposition \ref{Proposition3} are valid not only for $\lambda $- and $\mu $%
-based tests but also for any test invariant under left orthogonal
transformations of the observations ($X\mapsto QX,$ where $Q$ is a $p\times p
$ orthogonal matrix), and for any test invariant under multiplication by any
non-zero constant and left orthogonal transformations of the observations ($%
X\mapsto aQX,$ where $a\in \mathbb{R}_{0}^{+}$ and $Q$ is a $p\times p$
orthogonal matrix), respectively. Let $\left\Vert A\right\Vert _{F}=\mathrm{%
tr}\left( A^{\prime }A\right) $ and $\left\Vert A\right\Vert _{2}=\lambda
_{1}^{1/2}\left( A^{\prime }A\right) $ denote the Frobenius norm and the
spectral norm, respectively, of a matrix $A$. Let $H_{0}$ be the null
hypothesis $h_{1}=...=h_{r}=0,$ and let $H_{1}$ be any of the following
alternatives: $H_{1}:$ $h_{j}>0$ for some $j=1,...,r,$ or $H_{1}:\Sigma \neq
\sigma ^{2}I_{p},$ or $H_{1}:\left\{ \Sigma :\left\Vert \Sigma -\sigma
^{2}I_{p}\right\Vert _{F}>\varepsilon _{n,p}\right\} ,$ or $H_{1}:\left\{
\Sigma :\left\Vert \Sigma -\sigma ^{2}I_{p}\right\Vert _{2}>\varepsilon
_{n,p}\right\} $, where $\varepsilon _{n,p}$ is a positive constant that may
depend on $n$ and $p$.

\begin{proposition}
\label{Proposition4} For specified $\sigma ^{2}$, consider tests of $H_{0}$
against $H_{1}$ that are invariant with respect to the left orthogonal
transformations of the data $X=\left[ X_{1},...,X_{n}\right] .$ For any such
test, there exists a test based on $\lambda $ with the same power function.
Similarly, for unspecified $\sigma ^{2},$ consider tests that, in addition,
are invariant with respect to multiplication of the data $X$ by non-zero
constants. For any such test, there exists a test based on $\mu $ with the
same power function.
\end{proposition}

Examples of the former tests include the tests of $H_{0}:\Sigma =I$ studied
in Chen et al (2010) and Cai and Ma (2012). An example of the latter test is
the test of sphericity studied in Chen et al (2010). The tests studied in
Chen et al (2010) and Cai and Ma (2012) are invariant, although they are not 
$\lambda $- or $\mu $-based.

For $r=1,$ OMH show that the asymptotic power envelopes are closely
approached by the asymptotic powers of the $\lambda $- and $\mu $-based
likelihood ratio tests. Our goal here is to explore the asymptotic power of
those likelihood ratio tests for $r>1$. Unfortunately, as $r$ grows, it
becomes increasingly difficult to compute the asymptotic critical values for
the likelihood ratio tests by simulation. For example, $r=2$ requires
simulating a 2-dimensional Gaussian random field with the covariance
function and the mean function described in Theorem \ref{Theorem1}.

\begin{figure}[t!]
\centering
\includegraphics[width=4in]{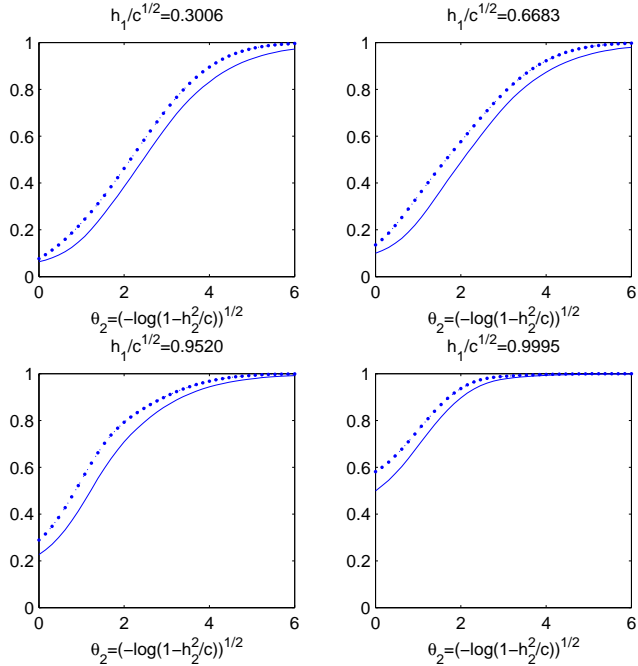}
\caption{Profiles of the asymptotic
power of the $\protect\lambda $-based LR test (solid lines) relative to the
asymptotic power envelope (dotted lines) for different values of $h_{1}/%
\protect\sqrt{c}$ under the alternative; $\protect\alpha =0.05.$ }
\label{powerlam_section}
\end{figure}

\begin{figure}[t!]
\centering
\includegraphics[width=4in]{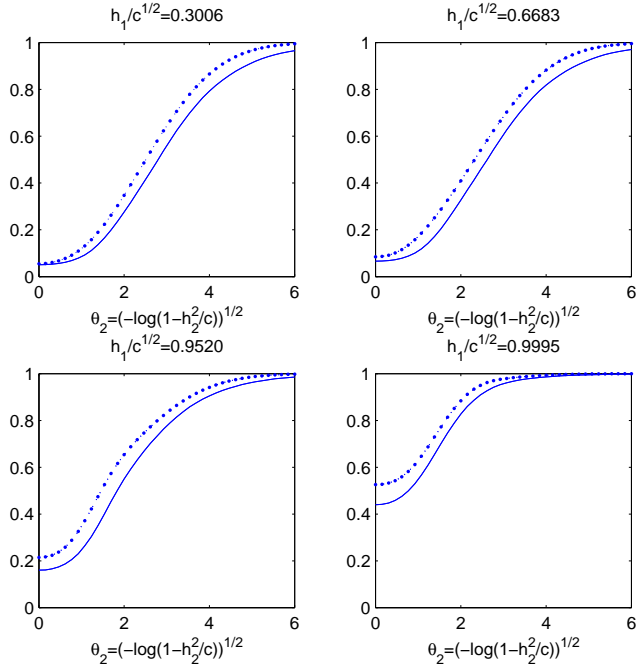}
\caption{Profiles of the asymptotic power
of the $\protect\mu $-based LR test (solid lines) relative to the asymptotic
power envelope (dotted lines) for different values of $h_{1}/\protect\sqrt{c}
$ under the alternative; $\protect\alpha =0.05.$ }
\label{powermu_section}
\end{figure}

For $r=2,$ Figure \ref{powerlam_section} shows sections of the power
envelope (dotted lines) and the power of the likelihood ratio test based on $%
\lambda $ for various fixed values of $h_{1}/\sqrt{c}$ under the
alternative. Figure \ref{powermu_section} shows the same plots for the tests
based on $\mu .$ To enhance readability, we use a different parametrization: 
$\theta _{j}=\sqrt{-\ln \left( 1-h_{j}^{2}/c\right) },$ $i=1,...,r$. As $%
h_{j}$ varies in the region of contiguity $\left[ 0,\sqrt{c}\right) ,$ $%
\theta _{j}$ spans the entire half-line $\left[ 0,\infty \right) .$ Note
that the asymptotic mean and autocovariance functions of the log likelihood
ratios derived in Theorem \ref{Theorem1} depend on $h_{j}$ only through $%
h_{j}/\sqrt{c}=\sqrt{1-e^{-\theta _{j}^{2}}}.$ Therefore, under the new
parametrization, they depend only on $\theta =\left( \theta _{1},...,\theta
_{r}\right) $. The parameter $\theta $ plays the classical role of a
\textquotedblleft local parameter\textquotedblright\ in our setting.


\begin{figure}[t!]
\centering
\includegraphics[width=4in]{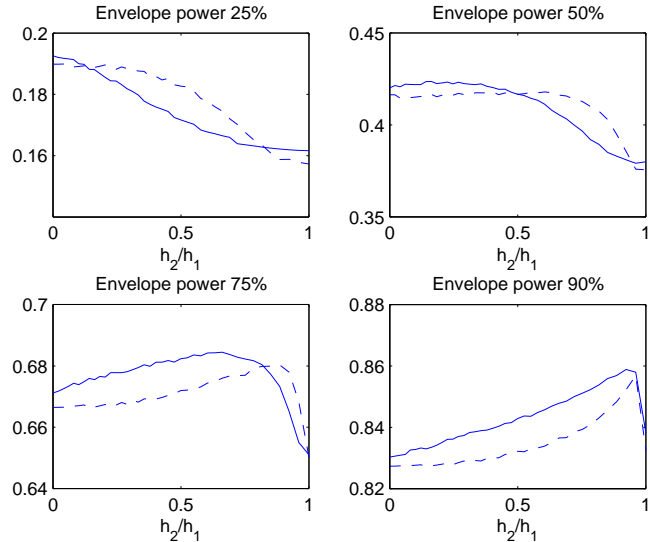}
\caption{Power of $\protect\lambda $-based
(solid lines) and $\protect\mu $-based (dashed lines) LR tests plotted
against $h_{2}/h_{1},$ where $\left( h_{1},h_{2}\right) $ are such that the
respective asymptotic power envelopes $\protect\beta _{\protect\lambda }(h)$
and $\protect\beta _{\protect\mu }(h)$ equal 25, 50, 75 and 90\%.}
\label{isoplot}
\end{figure}

Figure \ref{isoplot} further explores the relationship between the
asymptotic powers of the $\lambda $- and $\mu $-based LR test and the
corresponding asymptotic power envelopes when $r=2$. We pick all values of $%
h=\left( h_{1},h_{2}\right) $ satisfying inequality $h_{1}\geq h_{2}$ and
such that the asymptotic power envelope for $\lambda $-based tests is
exactly 25, 50, 75, and 90\%. Then, we compute and plot the corresponding
power of the $\lambda $-based LR test (solid lines) against $h_{2}/h_{1}.$
The dashed lines show similar graphs for the $\mu $-based LR test. The value 
$h_{2}/h_{1}=0$ corresponds to single-spiked alternatives $h_{1}>0,$~$%
h_{2}=0,$ the value $h_{2}/h_{1}=1$ corresponds to equi-spiked alternatives $%
h_{1}=h_{2}>0.$ The intermediate values of $h_{2}/h_{1}$ link the two
extreme cases. We do not consider values $h_{2}/h_{1}>1$, as the power
function is symmetric about the 45-degree line in the $\left(
h_{1},h_{2}\right) $ space.


Somewhat surprisingly, the power of the LR test along the set of
alternatives~$\left( h_{1},h_{2}\right) $ corresponding to the same values
of the asymptotic power envelope is not a monotone function of $h_{2}/h_{1}.$
Equi-spiked alternatives typically seem to be particularly difficult to
detect by the LR tests. However, for the set of alternatives corresponding
to an asymptotic power envelope value of 90\%, the single-spiked
alternatives are even harder to detect.

\begin{figure}[t!]
\centering
\includegraphics[width=4in]{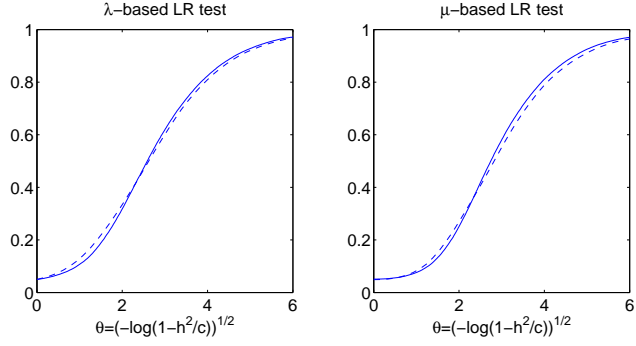}
\caption{Asymptotic power of the $%
\protect\lambda $-based (left panel) and $\protect\mu $-based (right panel)
LR tests. Solid line: power when $r=1$ is correctly assumed. Dashed line:
power when $r=2$ is incorrectly assumed.}
\label{wrongr}
\end{figure}

A natural question is: how does the asymptotic power of the $\lambda $- and $%
\mu $-based LR tests depend on the choice of $r$, that is, how do those
tests perform when the actual $r$ does not coincide with the value the test
statistic is based on? For example, to detect a single signal, one can, in
principle, use LR tests of the null hypothesis against alternatives with $%
r=1,r=2,$ etc. How does the asymptotic powers of such tests compare? Figure %
\ref{wrongr} reports the asymptotic powers of the $\lambda $- and $\mu $%
-based LR tests designed to detect alternatives with $r=1$ (solid line) and $%
r=2$ (dashed line), under single-spiked alternatives. As in Figures 2 and 3,
we use the parametrization $\theta =\sqrt{-\ln \left( 1-h^{2}/c\right) }$
for the single-spiked alternative. It appears that the two asymptotic powers
are very close to each other; interestingly, neither of them dominates the
other. Using LR tests designed against alternatives with $r>1$ seems to be
beneficial for detecting a single-spiked alternative with relatively small~$%
\theta $ (and $h$).


In the remaining part of this section, we consider examples of some of the
tests that have been proposed previously in the literature, and, in
Proposition \ref{Proposition5}, derive their asymptotic power functions.

\begin{example}[John's (1971) test of sphericity $H_{0}:\Sigma =\protect%
\sigma ^{2}I.$]
John (1971) proposes testing the sphericity hypothesis $\theta =0$ against
general alternatives based on the test statistic 
\begin{equation}
U=\frac{1}{p}\limfunc{tr}\left[ \left( \frac{\hat{\Sigma}}{\left( 1/p\right) 
\limfunc{tr}\left( \hat{\Sigma}\right) }-I_{p}\right) ^{2}\right] ,
\label{John's stat}
\end{equation}%
where $\hat{\Sigma}$ is the sample covariance matrix. He shows that, when $%
n>p,$ such a test is locally most powerful invariant. Ledoit and Wolf (2002)
study John's test when $n,p\rightarrow _{c}\infty $. They prove that, under
the null, $nU-p\overset{d}{\rightarrow }N\left( 1,4\right) .$ Hence, the
test with asymptotic size $\alpha $ rejects the null of sphericity whenever $%
\frac{1}{2}\left( nU-p-1\right) >\Phi ^{-1}\left( 1-\alpha \right) $.
\end{example}

\begin{example}[The Ledoit-Wolf (2002) test of $H_{0}:\Sigma =I.$]
Ledoit and Wolf (2002) propose 
\begin{equation}
W=\frac{1}{p}\limfunc{tr}\left[ \left( \hat{\Sigma}-I\right) ^{2}\right] -%
\frac{p}{n}\left[ \frac{1}{p}tr\hat{\Sigma}\right] ^{2}+\frac{p}{n}
\label{LW stat}
\end{equation}%
as a test statistic for testing the hypothesis that the population
covariance matrix is the unit matrix. They show that, under the null, $nW-p%
\overset{d}{\rightarrow }N\left( 1,4\right) .$ As in the previous example,
the null is rejected at asymptotic size $\alpha $ whenever $\frac{1}{2}%
\left( nW-p-1\right) >\Phi ^{-1}\left( 1-\alpha \right) .$
\end{example}

\begin{example}[The Bai et al. (2009) \textquotedblleft
corrected\textquotedblright\ LRT of $H_{0}:\Sigma =I$.]
When $n>p,$ Bai et al. (2009) propose to use a corrected version\vspace{%
-0.2in}%
\begin{equation*}
CLR=\limfunc{tr}\hat{\Sigma}-\ln \det \hat{\Sigma}-p-p\left( 1-\left( 1-%
\frac{n}{p}\right) \ln \left( 1-\frac{p}{n}\right) \right) \vspace{-0.2in}
\end{equation*}%
of the likelihood ratio statistic to test the equality of the population
covariance matrix to the identity matrix against general alternatives. Under
the null, $CLR\overset{d}{\rightarrow }N\left( -\frac{1}{2}\ln \left(
1-c\right) ,-2\ln \left( 1-c\right) -2c\right) $ (still, as $n,p\rightarrow
_{c}\infty $).\ The null hypothesis is rejected at asymptotic level $\alpha $
whenever $CLR+\frac{1}{2}\ln \left( 1-c\right) $ is larger than\linebreak\ $%
\left( -2\ln \left( 1-c\right) -2c\right) ^{1/2}\Phi ^{-1}\left( 1-\alpha
\right) $.
\end{example}

\begin{example}[Tracy-Widom-type tests of $H_{0}:\Sigma =I$.]
Let $\varphi \left( \lambda _{1},...,\lambda _{r}\right) $ be any function
of the $r$ largest eigenvalues increasing in all its arguments. The
asymptotic distribution of $\varphi \left( \lambda _{1},...,\lambda
_{r}\right) $ under the null is determined by the functional form of $%
\varphi \left( \cdot \right) $ and the fact that\vspace{-0.2in}%
\begin{equation}
\left( \sigma _{n,c}\left( \lambda _{1}-\nu _{c}\right) ,...,\sigma
_{n,c}\left( \lambda _{r}-\nu _{c}\right) \right) \overset{d}{\rightarrow }%
TW\left( r\right) ,\vspace{-0.2in}  \label{TW}
\end{equation}%
where TW$\left( r\right) $ denotes the $r$-dimensional Tracy-Widom law of
the first kind, $\sigma _{n,c}=n^{2/3}c^{1/6}\left( 1+\sqrt{c}\right)
^{-4/3} $ and $\nu _{c}=\left( 1+\sqrt{c}\right) ^{2}$. Call
Tracy-Widom-type tests all tests that reject the null whenever $\varphi
\left( \lambda _{1},...,\lambda _{r}\right) $ is larger than the
corresponding asymptotic critical value obtained from (\ref{TW}).
\end{example}

\begin{example}[The Cai-Ma (2012) minimax test of $H_{0}:\Sigma =I$.]
Cai and Ma (2012) propose to use a U-statistic\vspace{-0.2in}%
\begin{equation*}
T_{n}=\frac{2}{n\left( n-1\right) }\sum_{1\leq i<j\leq n}\ell \left(
X_{i},X_{j}\right) ,\vspace{-0.2in}
\end{equation*}%
where $\ell \left( X_{1},X_{2}\right) =\left( X_{1}^{\prime }X_{2}\right)
^{2}-\left( X_{1}^{\prime }X_{1}+X_{2}^{\prime }X_{2}\right) +p,$ to test
the hypothesis that the population covariance matrix is the unit matrix.
Under the null, as $n,p\rightarrow _{c}\infty $, $T_{n}\overset{d}{%
\rightarrow }N\left( 0,4c^{2}\right) .$ The null hypothesis is rejected at
asymptotic level $\alpha $ whenever $T_{n}$ is larger than $2\sqrt{p\left(
p+1\right) /n\left( n-1\right) }\Phi ^{-1}\left( 1-\alpha \right) $. Cai and
Ma (2012) show that this test is rate-optimal against general alternatives
from a minimax point of view.
\end{example}

Consider the tests described in Examples 1, 2, 3, 4 and 5, and denote by $%
\beta _{J}\left( h\right) ,$ $\beta _{LW}\left( h\right) ,$ $\beta
_{CLR}\left( h\right) ,$ $\beta _{CM}\left( h\right) $ and $\beta
_{TW}\left( h\right) $ their respective asymptotic powers at asymptotic
level $\alpha .$

\begin{proposition}
\label{Proposition5}\textit{The asymptotic power functions of the tests
described in Examples~1-5 are}%
\begin{eqnarray}
\beta _{TW}\left( h\right) \! &=&\!\alpha ,  \label{beta_TW} \\
\beta _{J}\left( h\right) \! &=&\!\beta _{LW}\left( h\right) \!=\!\beta
_{CM}\left( h\right) \!=\!1\!-\!\Phi \left( \Phi ^{-1}\left( 1\!-\!\alpha
\right) -\frac{1}{2}\sum\limits_{j=1}^{r}\frac{h_{j}^{2}}{c}\right) ,\text{
and}  \label{beta_J_LW} \\
\beta _{CLR}\left( h\right) \! &=&\!1\!-\!\Phi \left( \Phi ^{-1}\left(
1-\alpha \right) -\sum_{j=1}^{r}\frac{h_{j}-\ln \left( 1+h_{j}\right) }{%
\sqrt{-2\ln \left( 1-c\right) -2c}}\right) ,  \label{beta_CLR}
\end{eqnarray}%
for any $h=\left( h_{1},...,h_{r}\right) \neq 0$ such that $h_{j}\in \left[
0,\sqrt{c}\right) $ for $j=1,...,r.$
\end{proposition}

\begin{figure}[t!]
\centering
\includegraphics[width=4in]{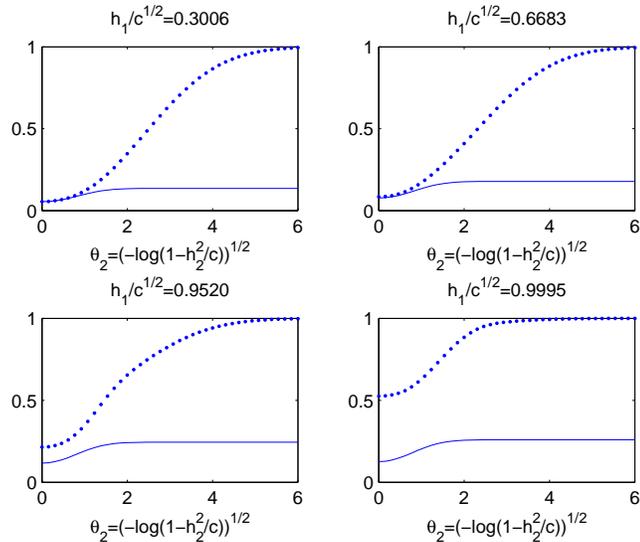}
\caption{Profiles of the asymptotic power
of John's test (solid lines) relative to the asymptotic power envelope $%
\protect\beta _{\protect\mu }$ (dotted lines) for different values of $h_{1}/%
\protect\sqrt{c}$ under the alternative; $\protect\alpha =0.05.$ }
\label{johnmu_section}
\end{figure}

\begin{figure}[t!]
\centering
\includegraphics[width=4in]{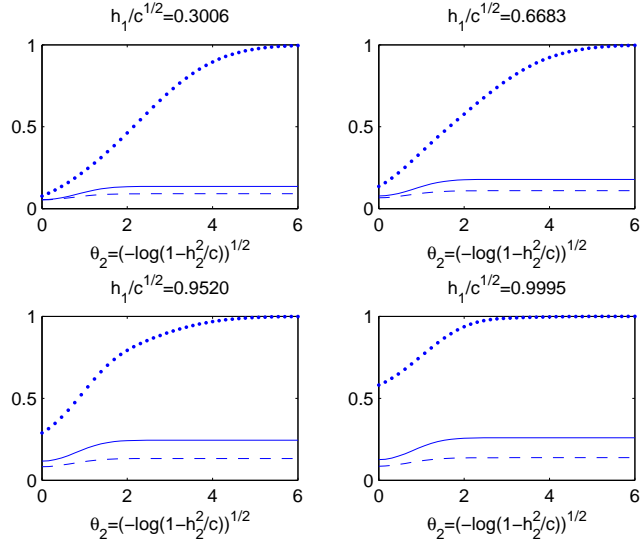}
\caption{Profiles of the asymptotic power
of the Ledoit-Wolf and Cai-Ma tests (solid lines) and the CLR test (dashed
lines, for $c=0.5$) relative to the asymptotic power envelope $\protect\beta %
_{\protect\lambda }$ (dotted lines) for different values of $h_{1}/\protect%
\sqrt{c}$ under the alternative; $\protect\alpha =0.05.$ }
\label{johnlam_section}
\end{figure}

Formula (\ref{beta_J_LW}) for $\beta _{CM}\left( h\right) $ directly follows
from Proposition 2 of Cai and Ma (2012). The proof of the other formulae
follows along the same lines as in the proof of Proposition 10 in OMH, and
is omitted. Except for the Tracy-Widom tests of Example 4, all those
asymptotic power functions are non-trivial. Figures \ref{johnmu_section} and~%
\ref{johnlam_section} compare these power functions to the corresponding
power envelopes for $r=2$. Since John's test is invariant with respect to
orthogonal transformations and scalings of the data, Figure \ref%
{johnmu_section} compares $\beta _{J}\left( h\right) $ (solid line) to the
power envelope $\beta _{\mu }\left( h\right) $ (dotted line). The
Ledoit-Wolf test, the \textquotedblleft corrected\textquotedblright\
likelihood ratio test, and the Cai-Ma test are invariant with respect to
orthogonal transformations of the data only, and Figure \ref{johnlam_section}
thus compares the asymptotic power functions $\beta _{LW}\left( h\right)
=\beta _{CM}\left( h\right) $ and $\beta _{CLR}\left( h\right) $ (solid and
dashed lines, respectively) to the power envelope $\beta _{\lambda }\left(
h\right) $ (dotted line). Note that $\beta _{CLR}\left( h\right) $ depends
on $c$. As $c$ converges to one, $\beta _{CLR}\left( h\right) $ converges to 
$\alpha ,$ which corresponds to the case of trivial power. As $c$ converges
to zero, $\beta _{CLR}\left( h\right) $ converges to $\beta _{LW}\left(
h\right) =\beta _{CM}\left( h\right) $. In Figure \ref{johnlam_section}, we
provide plots of~$\beta _{CLR}\left( h\right) $ that correspond to $c=0.5.$

These comparisons show that, contrary to our LR tests (see Figures \ref%
{powerlam_section} and \ref{powermu_section}), \textit{all} those tests
either have trivial power $\alpha $ (the Tracy-Widom ones), or power
functions that increase very slowly with $h_{1}$ and $h_{2},$ and lie very
far below the corresponding power envelope.


\section{Conclusion}

This paper extends Onatski, Moreira and Hallin's (2011) (OMH) study of the
power of high-dimensional sphericity tests to the case of multi-spiked
alternatives. We derive the asymptotic distribution of the log-likelihood
ratio process and use it to obtain simple analytical expressions for the
maximal asymptotic power envelope and for the asymptotic powers of several
tests proposed in the literature. These asymptotic powers turn out to be
very substantially below the envelope. We propose the likelihood ratio test
based on the data reduced to the eigenvalues of the sample covariance
matrix. Our computations show that the asymptotic power of this test is
close to the envelope.

\section{Appendix}

All convergence statements made below refer to the situation when $%
n_{p},p\rightarrow _{c}\infty $. We start with two auxiliary results.

\begin{lemma}
Let $d\left( \mu ,\nu \right) $ be the Dudley distance between measures $\mu 
$ and $\nu $ defined over $\left( \mathbb{R},\mathcal{B}\right) $:\vspace{%
-0.2in}%
\begin{equation*}
d\left( \mu ,\nu \right) =\sup \left\{ \left\vert \int f\left( d\mu -d\nu
\right) \right\vert :f\left( x\right) \leq 1\text{ and }\left\vert \frac{%
f\left( x\right) -f\left( y\right) }{x-y}\right\vert \leq 1,\forall x\neq
y\right\} .\vspace{-0.2in}
\end{equation*}%
There exists a constant $\tau >0$ such that $d\left( \hat{F}_{p}^{\lambda
},F_{p}^{MP}\right) =o\left( p^{-1}\log ^{\tau }p\right) $ a.s..
\end{lemma}

\textbf{Proof:} Let us denote the cumulative distribution function
corresponding to a measure $\mu $ as $F_{\mu }\left( x\right) .$ Further,
let us denote $\inf \left\{ \left\vert x_{2}-x_{1}\right\vert :\limfunc{supp}%
\left( \mu \right) \subseteq \left[ x_{1},x_{2}\right] \right\} $ as $%
\limfunc{diam}\left( \mu \right) .$ Consider the following three distances
between measures $\mu $ and $\nu :$ the Kolmogorov distance $k\left( \mu
,\nu \right) =\sup_{x}\left\vert F_{\mu }\left( x\right) -F_{\nu }\left(
x\right) \right\vert ,$ the Wasserstein distance $w\left( \mu ,\nu \right)
=\sup \left\{ \left\vert \int f\left( d\mu -d\nu \right) \right\vert
:\left\vert \frac{f\left( x\right) -f\left( y\right) }{x-y}\right\vert \leq
1,\forall x\neq y\right\} ,$ and the Kantorovich distance $\gamma \left( \mu
,\nu \right) =\int \left\vert F_{\mu }\left( x\right) -F_{\nu }\left(
x\right) \right\vert dx.$ As is well known (see, for example, exercise 1 on
p.425 of Dudley (2002)), $w\left( \mu ,\nu \right) =\gamma \left( \mu ,\nu
\right) .$ Therefore, we have$\vspace{-0.2in}$%
\begin{equation*}
d\left( \hat{F}_{p}^{\lambda },F_{p}^{MP}\!\right) \leq w\left( \hat{F}%
_{p}^{\lambda },F_{p}^{MP}\!\right) =\gamma \left( \hat{F}_{p}^{\lambda
},F_{p}^{MP}\!\right) \leq k\left( \hat{F}_{p}^{\lambda
},F_{p}^{MP}\!\right) \left( \limfunc{diam}\!\left( \hat{F}_{p}^{\lambda
}\right) +\limfunc{diam}\!\left( F_{p}^{MP}\right) \right) .\vspace{-0.2in}
\end{equation*}%
As follows from Theorem 1.1 of G\"{o}tze and Tikhomirov (2011), there exists
a constant $\tau >0$ such that $\sum_{p=1}^{\infty }\Pr \left( k\left( \hat{F%
}_{p}^{\lambda },F_{p}^{MP}\right) >\varepsilon p^{-1}\log ^{\tau }p\right)
<\infty $ for all $\varepsilon >0$. Thus, $k\left( \hat{F}_{p}^{\lambda
},F_{p}^{MP}\right) =o\left( p^{-1}\log ^{\tau }p\right) $ a.s.. Since $%
\limfunc{diam}\left( F_{p}^{MP}\right) $ is $O(1)$ and $\limfunc{diam}\left( 
\hat{F}_{p}^{\lambda }\right) -\limfunc{diam}\left( F_{p}^{MP}\right)
\rightarrow 0$ a.s., the result follows.$\square $

\begin{corollary}
Suppose that a sequence of functions $\left\{ f_{p}(\lambda )\right\} $ is
bounded Lipshitz on $\limfunc{supp}\left( F_{p}^{MP}\right) \cup \limfunc{%
supp}\left( \hat{F}_{p}^{\lambda }\right) $, uniformly over all sufficiently
large $p,$ a.s.. Then $\left\vert \int f_{p}(\lambda )\mathrm{d}\left( \hat{F%
}_{p}^{\lambda }\left( \lambda \right) -F_{p}^{MP}\left( \lambda \right)
\right) \right\vert =o\left( p^{-1/2}\right) ,$ a.s..
\end{corollary}

\subsection{Proof of Proposition \protect\ref{Proposition2}}

Let us denote the integral $\int\limits_{\mathcal{O}\left( p\right) }e^{p%
\limfunc{tr}\left( \Theta _{p}Q^{\prime }\Lambda _{p}Q\right) }\left( 
\mathrm{d}Q\right) $ as $I_{p}\left( \Theta _{p},\Lambda _{p}\right) $. As
explained in Guionnet and Maida (2005, p.454), we can write$\vspace{-0.2in}$%
\begin{equation}
I_{p}\left( \Theta _{p},\Lambda _{p}\right) =\mathbb{E}_{\Lambda _{p}}\exp
\left\{ p\sum_{j=1}^{r}\theta _{pj}\frac{\tilde{g}^{(j)\prime }\Lambda _{p}%
\tilde{g}^{(j)}}{\tilde{g}^{(j)\prime }\tilde{g}^{(j)}}\right\} ,\vspace{%
-0.2in}  \label{First Ip}
\end{equation}%
where $\mathbb{E}_{\Lambda _{p}}$ denotes the expectation conditional on $%
\Lambda _{p},$ and the $p$-dimensional vectors $\left( \tilde{g}^{(1)},...,%
\tilde{g}^{(r)}\right) $ are obtained from standard Gaussian $p$-dimensional
vectors $\left( g^{(1)},...,g^{\left( r\right) }\right) $, independent from $%
\Lambda _{p}$, by a Schmidt orthogonalization procedure. More precisely, we
have $\tilde{g}^{(j)}=\sum_{k=1}^{j}A_{jk}g^{(k)}$, where $A_{jj}=1$ and$%
\vspace{-0.2in}$%
\begin{equation}
\sum_{k=1}^{j-1}A_{jk}g^{(k)\prime }g^{(t)}=-g^{(j)\prime }g^{(t)}\text{ for 
}t=1,...,j-1\text{.}\vspace{-0.2in}  \label{A definition}
\end{equation}

In the spirit of the proof of Guionnet and Maida's (2005) Theorem 3, define$%
\vspace{-0.2in}$%
\begin{equation}
\gamma _{p1}^{(j,s)}=\!\sqrt{p}\left( \frac{1}{p}g^{(j)\prime
}g^{(s)}\!-\!\delta _{js}\right) \text{ and }\gamma _{p2}^{(j,s)}=\!\sqrt{p}%
\left( \frac{1}{p}g^{(j)\prime }\Lambda _{p}g^{(s)}\!-\!v_{pj}\delta
_{js}\right) ,\vspace{-0.2in}  \label{gamma definition}
\end{equation}%
where $\delta _{js}=\mathbf{1}\left\{ j=s\right\} $ stands for the classical
Kronecker symbol. As will be shown below, after an appropriate change of
measure, $\gamma _{p1}^{(j,s)}$ and $\gamma _{p2}^{(j,s)}$ are
asymptotically centered Gaussian. Expressing the exponent in (\ref{First Ip}%
) as a function of $\gamma _{p1}^{(j,s)}$ and~$\gamma _{p2}^{(j,s)},$
changing the measure of integration, and using the asymptotic Gaussianity
will establish the proposition.

Let $\gamma _{p}\!=\!\left( \!\gamma _{p}^{(1,1)}\!,...,\gamma
_{p}^{(r,1)}\!,\gamma _{p}^{(2,2)}\!,...,\gamma _{p}^{(r,2)}\!,\gamma
_{p}^{(3,3)}\!,...,\gamma _{p}^{(r,r)}\!\right) ^{\prime }$, where $\gamma
_{p}^{(j,s)}\!=\!\left( \gamma _{p1}^{(j,s)}\!,\gamma _{p2}^{(j,s)}\right) $%
. Using this notation, (\ref{First Ip}), (\ref{A definition}), and (\ref%
{gamma definition}), we get, after some algebra,$\vspace{-0.2in}$%
\begin{equation}
I_{p}\left( \Theta _{p},\Lambda _{p}\right) \!=\!\int \!f_{p,\theta }\left(
\gamma _{p}\right) e^{p\sum_{j=1}^{r}\theta _{pj}\left( v_{pj}+\hat{\gamma}%
_{p}^{(j,j)}-v_{pj}\gamma _{p}^{(j,j)}\right) }\prod_{j=1}^{r}\prod_{i=1}^{p}%
\mathrm{d}\mathbb{P}\left( g_{i}^{(j)}\right) ,\vspace{-0.2in}  \label{IMM1}
\end{equation}%
where $\mathbb{P}$ is the standard Gaussian probability measure, and$\vspace{%
-0.2in}$%
\begin{eqnarray}
f_{p,\theta }\left( \gamma _{p}\right)  &=&\exp \left\{ \sum_{j=1}^{r}\theta
_{pj}\frac{N_{1j}\!+...\!+\!N_{6j}}{D_{j}}\right\} \text{ with}  \label{fp}
\\
N_{1j} &=&-\gamma _{p1}^{(j,j)}\left( \gamma _{p2}^{(j,j)}-v_{pj}\gamma
_{p1}^{(j,j)}\right) ,  \notag \\
N_{2j} &=&\gamma _{p1}^{\left( j,1:j-1\right) \prime }\left(
G_{p1}^{(j)}+I\right) ^{-1}\left( G_{p2}^{(j)}+W_{pj}\right) \left(
G_{p1}^{(j)}+I\right) ^{-1}\gamma _{p1}^{\left( j,1:j-1\right) },  \notag \\
N_{3j} &=&-2\gamma _{p1}^{\left( j,1:j-1\right) \prime }\left(
G_{p1}^{(j)}+I\right) ^{-1}\gamma _{p2}^{\left( j,1:j-1\right) },  \notag \\
N_{4j} &=&v_{pj}\gamma _{p1}^{\left( j,1:j-1\right) \prime }\left(
G_{p1}^{(j)}+I\right) ^{-1}\gamma _{p1}^{\left( j,1:j-1\right) },  \notag \\
N_{5j} &=&p^{-1/2}\gamma _{p2}^{(j,j)}\gamma _{p1}^{\left( j,1:j-1\right)
\prime }\left( G_{p1}^{(j)}+I\right) ^{-1}\gamma _{p1}^{\left(
j,1:j-1\right) },  \notag \\
N_{6j} &=&-p^{-1/2}v_{pj}\gamma _{p1}^{\left( 1:j-1,j\right) \prime }\left(
G_{p1}^{(j)}+I\right) ^{-1}\gamma _{p1}^{\left( 1:j-1,j\right) }\gamma
_{p1}^{(j,j)},\text{ and}  \notag \\
D_{j} &=&1+p^{-1/2}\gamma _{p1}^{(j,j)}-p^{-1}\gamma _{p1}^{\left(
j,1:j-1\right) \prime }\left( G_{p1}^{(j)}+I\right) ^{-1}\gamma
_{p1}^{\left( j,1:j-1\right) }\text{,}\vspace{-0.3in}  \notag
\end{eqnarray}%
where $G_{pi}^{(j)}$ is a $\left( j-1\right) \times \left( j-1\right) $
matrix with $\left( k,s\right) $-th element $p^{-1/2}\gamma _{pi}^{(k,s)},$%
\vspace{-0.2in}%
\begin{equation*}
W_{pj}=\mathrm{diag}\left( v_{p1},...,v_{p,j-1}\right) ,\text{ and }\gamma
_{pi}^{\left( j,1:j-1\right) }=\left( \gamma _{pi}^{\left( j,1\right)
},...,\gamma _{pi}^{\left( j,j-1\right) }\right) ^{\prime }.\vspace{-0.2in}
\end{equation*}

Next, define the event$\vspace{-0.2in}$%
\begin{equation*}
B_{M,M^{\prime }}=\left\{ \left\vert \gamma _{p1}^{(j,s)}\right\vert \leq M%
\text{ and }\left\vert \gamma _{p2}^{(j,s)}\right\vert \leq M^{\prime }\text{
for all }j,s=1,...,r\right\} ,\vspace{-0.2in}
\end{equation*}%
where $M$ and $M^{\prime }$ are positive parameters to be specified later.
Somewhat abusing notation, we will also refer to $B_{M,M^{\prime }}$ as a
rectangular region in $R^{r^{2}+r}$ that consists of vectors with odd
coordinates in $\left( -M,M\right) $ and even coordinates in $\left(
-M^{\prime },M^{\prime }\right) $. Let$\vspace{-0.2in}$%
\begin{equation*}
I_{p}^{M,M^{\prime }}\left( \Theta _{p},\Lambda _{p}\right) \!=\!\int \!%
\mathbf{1}\left\{ \!B_{M,M^{\prime }}\!\right\} f_{p,\theta }\left( \gamma
_{p}\right) e^{p\sum_{j=1}^{r}\theta _{pj}\left( v_{pj}+\hat{\gamma}%
_{p}^{(j,j)}-v_{pj}\gamma _{p}^{(j,j)}\right) }\prod_{j=1}^{r}\prod_{i=1}^{p}%
\mathrm{d}\mathbb{P}\left( g_{i}^{(j)}\right) ,\vspace{-0.2in}
\end{equation*}%
where $\mathbf{1}\left\{ \cdot \right\} $ denotes the indicator function.
Below, we establish the asymptotic behavior of $I_{p}^{M,M^{\prime }}\left(
\Theta _{p},\Lambda _{p}\right) $ as first $p,$ and then $M$ and $M^{\prime }
$, diverge to infinity. We then show that the asymptotics of $%
I_{p}^{M,M^{\prime }}\left( \Theta _{p},\Lambda _{p}\right) $ and $%
I_{p}\left( \Theta _{p},\Lambda _{p}\right) $ coincide.

Consider infinite arrays $\left\{ \mathbb{P}_{pi}^{(j)},p=1,2,...;i=1,...,p%
\right\} ,$ $j=1,...,r,$ of random centered Gaussian measures$\vspace{-0.2in}
$%
\begin{equation*}
\mathrm{d}\mathbb{P}_{pi}^{(j)}\left( x\right) =\sqrt{\frac{1+2\theta
_{pj}v_{pj}-2\theta _{pj}\lambda _{pi}}{2\pi }}e^{-\frac{1}{2}\left(
1+2\theta _{pj}v_{pj}-2\theta _{pj}\lambda _{pi}\right) x^{2}}\mathrm{d}x.%
\vspace{-0.2in}
\end{equation*}%
Since $v_{pj}=R_{p}^{MP}\left( 2\theta _{pj}\right) =1/\left( 1-2\theta
_{pj}c_{p}\right) $ and $2\theta _{pj}\in \Omega _{\varepsilon \eta }$,
there exists $\hat{\varepsilon}>0$ such that, for sufficiently large $p,%
\vspace{-0.2in}$%
\begin{eqnarray*}
v_{pj}+\frac{1}{2\theta _{pj}} &>&\left( 1+\sqrt{c}\right) ^{2}+\hat{%
\varepsilon}\text{ when }\theta _{pj}>0\text{ and} \\
v_{pj}+\frac{1}{2\theta _{pj}} &<&-\hat{\varepsilon}\text{ when }\theta
_{pj}<0.\vspace{-0.2in}
\end{eqnarray*}%
Recall that $\lambda _{pp}\geq 0,$ and $\lambda _{p1}\rightarrow \left( 1+%
\sqrt{c}\right) ^{2}$ a.s.. Therefore, still a.s., for sufficiently large $%
p, $ $v_{pj}+\frac{1}{2\theta _{pj}}>\lambda _{p1}$ when~$\theta _{pj}>0$
and $v_{pj}+\frac{1}{2\theta _{pj}}<\lambda _{pp}$ when $\theta _{pj}<0.$
Hence, the measures $\mathbb{P}_{pi}^{(j)}$ are a.s. well defined for
sufficiently large $p$. Whenever $\mathbb{P}_{pi}^{(j)}$ is not well
defined, we re-define it arbitrarily.

We have$\vspace{-0.2in}$%
\begin{equation}
I_{p}^{M,M^{\prime }}\left( \Theta _{p},\Lambda _{p}\right)
=e^{p\sum_{j=1}^{r}\left[ \theta _{pj}v_{pj}-\frac{1}{2p}\sum_{i=1}^{p}\ln
\left( 1+2\theta _{pj}v_{pj}-2\theta _{pj}\lambda _{pi}\right) \right]
}J_{p}^{M,M^{\prime }},\vspace{-0.2in}  \label{IMMnow}
\end{equation}%
where$\vspace{-0.2in}$%
\begin{equation}
J_{p}^{M,M^{\prime }}=\int \!\mathbf{1}\left\{ B_{M,M^{\prime }}\right\}
f_{p,\theta }\left( \gamma _{p}\right) \prod_{j=1}^{r}\prod_{i=1}^{p}\mathrm{%
d}\mathbb{P}_{pi}^{(j)}\left( g_{i}^{(j)}\right) .  \label{JMM_first}
\end{equation}%
We now show that, under $\prod_{j=1}^{r}\prod_{i=1}^{p}d\mathbb{P}%
_{pi}^{(j)}\left( g_{i}^{(j)}\right) $, $\gamma _{p}$ a.s. converges in
distribution to a centered $r^{2}+r$-dimensional Gaussian vector, so that $%
J_{p}^{M,M^{\prime }}$ is asymptotically equivalent to an integral with
respect to a Gaussian measure on $\mathbb{R}^{r^{2}+r}.$

First, let us find the mean $\mathbb{E}_{p}\gamma _{p}$, and the variance $%
\mathbb{V}_{p}\gamma _{p}$ of $\gamma _{p}$ under measure $%
\prod_{j=1}^{r}\!\prod_{i=1}^{p}\!\mathrm{d}\mathbb{P}_{pi}^{(j)}\left(
\!g_{i}^{(j)}\!\right) $. Note that $\!\mathbb{V}_{p}\gamma _{p}\!=\!\mathrm{%
diag}\left( \!\mathbb{V}_{p}\gamma _{p}^{(1,1)}\!,\!\mathbb{V}_{p}\gamma
_{p}^{(2,1)}\!,\!...\!,\!\mathbb{V}_{p}\gamma _{p}^{(r,r)}\!\right) $ and $%
\mathbb{E}_{p}\gamma _{p}\!=\!\left( \!\mathbb{E}_{p}\gamma _{p}^{(1,1)}\!,%
\mathbb{E}_{p}\gamma _{p}^{(2,1)}\!,\!...,\mathbb{E}_{p}\gamma
_{p}^{(r,r)}\right) \!^{\prime }$. With probability one, for sufficiently
large $p,$ we have$\vspace{-0.2in}$%
\begin{eqnarray*}
\mathbb{E}_{p}\gamma _{p1}^{(k,s)} &=&\sqrt{p}\delta _{ks}\left( \frac{1}{p}%
\sum_{i=1}^{p}\frac{1}{\left( 1+2\theta _{pk}v_{pk}-2\theta _{pk}\lambda
_{pi}\right) }-1\right) \\
&=&\sqrt{p}\delta _{ks}\int \frac{\left( 2\theta _{pk}\right) ^{-1}}{%
K_{p}^{MP}\left( 2\theta _{pk}\right) -\lambda }\mathrm{d}\left( \hat{F}%
_{p}^{\lambda }\left( \lambda \right) -F_{p}^{MP}\left( \lambda \right)
\right) ,\vspace{-0.2in}
\end{eqnarray*}%
which, by Corollary 1, is $o\left( 1\right) $ uniformly in $2\theta _{pk}\in
\Omega _{\varepsilon \eta }$, a.s.. That Corollary 1 can be applied here
follows from the form of expression (\ref{K transform}) for $%
K_{p}^{MP}\left( x\right) $. Similarly,$\vspace{-0.2in}$%
\begin{equation*}
\mathbb{E}_{p}\gamma _{p2}^{(k,s)}\!=\!\sqrt{p}\frac{\delta _{ks}}{2\theta
_{pk}}\int \!\frac{K_{p}^{MP}\left( 2\theta _{pk}\right) }{K_{p}^{MP}\left(
2\theta _{pk}\right) -\lambda }\mathrm{d}\!\left( \!\hat{F}_{p}^{\lambda
}\left( \lambda \right) \!-\!F_{p}^{MP}\left( \lambda \right) \!\right)
\!=\!o\left( 1\right) \vspace{-0.2in}
\end{equation*}%
uniformly in $2\theta _{pk},2\theta _{ps}\in \Omega _{\varepsilon \eta }$,
a.s.. Thus,$\vspace{-0.2in}$%
\begin{equation}
\sup_{\left\{ 2\theta _{pj}\in \Omega _{\varepsilon \eta },j\leq r\right\} }%
\mathbb{E}_{p}\gamma _{p}=o\left( 1\right) \text{ a.s.}.\vspace{-0.2in}
\label{EgammaTil}
\end{equation}

Next, with probability one, for sufficiently large $p$ we have$\vspace{-0.2in%
}$%
\begin{equation*}
\mathbb{V}_{p}\gamma _{p1}^{(k,s)}=\frac{1}{p}\sum_{i=1}^{p}\frac{1+\delta
_{ks}}{\left( 1+2\theta _{pk}v_{pk}-2\theta _{pk}\lambda _{pi}\right) \left(
1+2\theta _{ps}v_{ps}-2\theta _{ps}\lambda _{pi}\right) }\text{.}\vspace{%
-0.2in}
\end{equation*}%
Let $\hat{H}_{p,ks}^{(2)}=\int \frac{\mathrm{d}\hat{F}_{p}^{\lambda }\left(
\lambda \right) }{\left( \!K_{p}^{MP}\!\left( \!2\theta _{pk}\!\right)
\!-\!\lambda \!\right) \left( \!K_{p}^{MP}\!\left( \!2\theta _{ps}\!\right)
\!-\!\lambda \!\right) }$ and $H_{p,ks}^{(2)}=\int \frac{\mathrm{d}%
F_{p}^{MP}\left( \lambda \right) }{\left( \!K_{p}^{MP}\!\left( \!2\theta
_{pk}\!\right) \!-\!\lambda \!\right) \left( \!K_{p}^{MP}\!\left( \!2\theta
_{ps}\!\right) \!-\!\lambda \!\right) }$. Then, using Corollary 1, we get$%
\vspace{-0.2in}$%
\begin{equation*}
\mathbb{V}_{p}\gamma _{p1}^{(k,s)}=\frac{1+\delta _{ks}}{4\theta _{pk}\theta
_{ps}}\hat{H}_{p,ks}^{(2)}=\frac{1+\delta _{ks}}{4\theta _{pk}\theta _{ps}}%
H_{p,ks}^{(2)}+o(1)\text{ a.s.},\vspace{-0.2in}
\end{equation*}%
uniformly in $2\theta _{pk},2\theta _{ps}\in \Omega _{\varepsilon \eta }$.
Similarly, we have$\vspace{-0.2in}$%
\begin{eqnarray*}
\mathbb{V}_{p}\gamma _{p2}^{(k,s)}\!\! &=&\!\!\frac{1}{p}\!\sum_{i=1}^{p}\!%
\frac{\lambda _{pi}^{2}\left( 1+\delta _{ks}\right) }{\left( \!1\!+\!2\theta
_{pk}v_{pk}\!-\!2\theta _{pk}\lambda _{pi}\!\right) \left( \!1\!+\!2\theta
_{ps}v_{ps}\!-\!2\theta _{ps}\lambda _{pi}\!\right) } \\
&=&\!\!\frac{1\!+\!\delta _{ks}}{4\theta _{pk}\theta _{ps}}\left(
\!1\!\!+\!\!K_{p}^{MP}\left( \!2\theta _{ps}\!\right) \!K_{p}^{MP}\left(
\!2\theta _{pk}\!\right) H_{p,ks}^{(2)}\!\!-\!\!2\theta
_{pk}K_{p}^{MP}\left( \!2\theta _{pk}\!\right) \!\!-\!\!2\theta
_{ps}K_{p}^{MP}\left( \!2\theta _{ps}\!\right) \!\right) \!+\!o(\!1\!),%
\vspace{-0.2in}
\end{eqnarray*}%
and$\vspace{-0.2in}$%
\begin{eqnarray*}
\mathbb{C}\mathrm{ov}_{p}\left( \!\gamma _{p1}^{(k,s)}\!,\gamma
_{p2}^{(k,s)}\!\right) \! &=&\!\frac{1}{p}\sum_{i=1}^{p}\frac{\lambda
_{pi}\left( 1+\delta _{ks}\right) }{\left( 1+2\theta _{pk}v_{pk}-2\theta
_{pk}\lambda _{pi}\right) \left( 1+2\theta _{ps}v_{ps}-2\theta _{ps}\lambda
_{pi}\right) } \\
&=&\frac{\left( 1+\delta _{ks}\right) }{4\theta _{pk}\theta _{ps}}\left(
K_{p}^{MP}\left( 2\theta _{ps}\right) H_{p,ks}^{(2)}\!-\!2\theta
_{pk}\right) +o(1),\vspace{-0.2in}
\end{eqnarray*}%
uniformly in $2\theta _{pk},2\theta _{ps}\in \Omega _{\varepsilon \eta }$,
a.s..

A straightforward calculation, using formula (\ref{K transform}), shows that$%
\vspace{-0.2in}$%
\begin{equation}
H_{p,ks}^{(2)}=\left( \frac{1}{4\theta _{pk}\theta _{ps}}-c_{p}v_{pk}v_{sk}%
\right) ^{-1}\text{, and }\mathbb{V}_{p}\gamma
_{p}^{(k,s)}=V_{p}^{(k,s)}+o(1),\vspace{-0.2in}  \label{VarGammaTil}
\end{equation}%
uniformly in $2\theta _{pk},2\theta _{ps}\in \Omega _{\varepsilon \eta }$,
a.s., where the matrix $V_{p}^{(k,s)}$ has elements$\vspace{-0.2in}$%
\begin{eqnarray}
V_{p,11}^{(k,s)}\! &=&\!\left( 1\!+\!\delta _{ks}\right) \left(
\!1\!-\!4\theta _{pk}v_{pk}\theta _{ps}v_{sk}c_{p}\!\right) ^{-1},
\label{V11} \\
V_{p,12}^{(k,s)}\! &=&\!V_{p,21}^{(k,s)}\!=\!\left( 1\!+\!\delta
_{ks}\right) v_{pk}v_{sk}\left( \!1\!-\!4\theta _{pk}v_{pk}\theta
_{ps}v_{sk}c_{p}\!\right) ^{-1},\text{ and}  \label{V12} \\
V_{p,22}^{(k,s)}\! &=&\!\left( 1\!+\!\delta _{ks}\right) \left[
\!c_{p}v_{pk}v_{sk}\!+\!v_{pk}^{2}v_{sk}^{2}\left( \!1\!-\!4\theta
_{pk}v_{pk}\theta _{ps}v_{sk}c_{p}\!\right) ^{-1}\right] \text{.}\vspace{%
-0.25in}  \label{V22}
\end{eqnarray}%
This implies that$\vspace{-0.2in}$%
\begin{equation}
\det \left( V_{p}^{(k,s)}\right) \!=\!\prod_{k\geq s}^{r}\left( 1\!+\!\delta
_{ks}\right) ^{2}c_{p}v_{pk}v_{sk}\left( \!1\!-\!4\theta _{pk}v_{pk}\theta
_{ps}v_{sk}c_{p}\!\right) ^{-1},\vspace{-0.2in}  \label{DetVariance}
\end{equation}%
which is bounded away from zero and infinity for sufficiently large $p$,
uniformly over $\left\{ 2\theta _{pj}\in \Omega _{\varepsilon \eta },j\leq
r\right\} $, a.s..

By construction, $\gamma _{p}$ is a sum of $p$ independent random vectors
having uniformly bounded third and fourth absolute moments under measure $%
\prod_{j=1}^{r}\prod_{i=1}^{p}\mathrm{d}\mathbb{P}_{pi}^{(j)}\left(
g_{i}^{(j)}\right) .$ Therefore, a central limit theorem applies. Moreover,
since the function $f_{p,\theta }\left( \gamma _{p}\right) $ is Lipshitz
over $B_{M,M^{\prime }},$ uniformly in $\left\{ 2\theta _{pj}\in \Omega
_{\varepsilon \eta },j\leq r\right\} ,$ Theorem 13.3 of Bhattacharya and Rao
(1976), which describes the accuracy of the Gaussian approximations to
integrals of the form (\ref{JMM_first}) in terms of the oscillation measures
of the integrand, implies that$\vspace{-0.2in}$%
\begin{equation}
J_{p}^{M,M^{\prime }}=\int_{B_{M,M^{\prime }}}f_{p,\theta }\left( x\right) 
\mathrm{d}\Phi \left( x;\mathbb{E}_{p}\gamma _{p},\mathbb{V}_{p}\gamma
_{p}\right) +o_{M,M^{\prime }}\left( 1\right) ,\vspace{-0.2in}
\label{JMM_second}
\end{equation}%
where $\Phi \left( x;\mathbb{E}_{p}\gamma _{p},\mathbb{V}_{p}\gamma
_{p}\right) $ denotes the Gaussian distribution function with mean $\mathbb{E%
}_{p}\gamma _{p}$ and variance $\mathbb{V}_{p}\gamma _{p},$ and $%
o_{M,M^{\prime }}\left( 1\right) $ converges to zero uniformly in $\left\{
2\theta _{pj}\in \Omega _{\varepsilon \eta },j\leq r\right\} $ as~$%
p\rightarrow \infty ,$ a.s.. The rate of such a convergence may depend on
the values of $M$ and $M^{\prime }.$

Note that, in $B_{M,M^{\prime }},$ as $p\rightarrow \infty ,$ the difference 
$f_{p,\theta }\left( \gamma _{p}\right) -\overline{f}_{p,\theta }\left(
\gamma _{p}\right) $ converges to zero uniformly over $\left\{ 2\theta
_{pj}\in \Omega _{\varepsilon \eta },j\leq r\right\} ,$ where$\vspace{-0.2in}
$%
\begin{eqnarray}
\overline{f}_{p,\theta }\left( \gamma _{p}\right)  &=&\exp \left\{
\sum_{j=1}^{r}\theta _{pj}\left( \bar{N}_{1j}\!+...\!+\!\bar{N}_{4j}\right)
\right\} ,\text{ with}  \label{fibar} \\
\bar{N}_{1j} &=&-\gamma _{1}^{(j,j)}\left( \gamma _{2}^{(j,j)}-v_{pj}\gamma
_{1}^{(j,j)}\right) ,\bar{N}_{2j}=\gamma _{1}^{\left( j,1:j-1\right) \prime
}W_{pj}\gamma _{1}^{\left( j,1:j-1\right) },  \notag \\
\bar{N}_{3j} &=&-2\gamma _{1}^{\left( j,1:j-1\right) \prime }\gamma
_{2}^{\left( j,1:j-1\right) },\text{ and }\bar{N}_{4j}=v_{pj}\gamma
_{1}^{\left( j,1:j-1\right) \prime }\gamma _{1}^{\left( j,1:j-1\right) }.%
\vspace{-0.3in}  \notag
\end{eqnarray}%
Such a convergence, together with (\ref{EgammaTil}), (\ref{VarGammaTil}),
and (\ref{JMM_second}) implies that$\vspace{-0.2in}$%
\begin{equation}
J_{p}^{M,M^{\prime }}=\int_{B_{M,M^{\prime }}}\overline{f}_{p,\theta }\left(
x\right) \mathrm{d}\Phi \left( x;0,V_{p}\right) +o_{M,M^{\prime }}\left(
1\right) ,\vspace{-0.2in}  \label{JMM_third}
\end{equation}%
where $V_{p}=\mathrm{diag}\left(
V_{p}^{(1,1)},V_{p}^{(2,1)},...,V_{p}^{(r,r)}\right) $.

Note that the difference $\int\limits_{B_{M,M^{\prime }}}\!\!\overline{f}%
_{p,\theta }\left( x\right) \mathrm{d}\Phi \left( \!x;\!0,\!V_{p}\!\right)
-\int\limits_{\mathbb{R}^{r^{2}+r}}\!\!\overline{f}_{p,\theta }\left(
x\right) \mathrm{d}\Phi \left( \!x;\!0\!,\!V_{p}\!\right) $ converges to
zero as $M,M^{\prime }\rightarrow \infty $, uniformly in $p$ for $p$
sufficiently large. On the other hand,$\vspace{-0.2in}$%
\begin{equation}
\int\limits_{\mathbb{R}^{r^{2}+r}}\!\!\overline{f}_{p,\theta }\left(
x\right) \mathrm{d}\Phi \left( x;0,V_{p}\right)
\!\!=\!\!\prod_{j=1}^{r}\prod_{s=1}^{j}\!\int\limits_{\mathbb{R}^{2}}\!\!%
\frac{\exp \left[ -\frac{1}{2}y^{\prime }\left( W_{p}^{(j,s)}\right) ^{-1}y%
\right] }{2\pi \!\sqrt{\det \left( V_{p}^{(j,s)}\right) }}\mathrm{d}y,%
\vspace{-0.2in}  \label{MMintegral}
\end{equation}%
where$\vspace{-0.2in}$%
\begin{equation*}
\left( W_{p}^{(j,s)}\right) ^{-1}=\left( V_{p}^{(j,s)}\right) ^{-1}+\left(
1+\delta _{js}\right) ^{-1}\left( 
\begin{array}{cc}
-2\theta _{pj}\left( v_{pj}+v_{ps}\right) & 2\theta _{pj} \\ 
2\theta _{pj} & 0%
\end{array}%
\right) \text{.}\vspace{-0.2in}
\end{equation*}%
Using (\ref{V11}-\ref{V22}), we verify that, for sufficiently large $p,$ $%
W_{p}^{(j,s)}$ is a.s. positive definite, and$\vspace{-0.2in}$%
\begin{eqnarray}
\det \left( W_{p}^{(j,s)}\right) \! &=&\!\left( 1+\delta _{js}\right)
^{2}c_{p}v_{pj}v_{ps},\text{ and}  \label{determinants} \\
\det \left( V_{p}^{(j,s)}\right) \! &=&\!\left( 1+\delta _{js}\right)
^{2}c_{p}v_{pj}v_{ps}\left( 1\!-\!4\left( \!\theta _{pj}v_{pj}\!\right)
\left( \!\theta _{ps}v_{ps}\!\right) c_{p}\right) ^{-1}.\vspace{-0.25in}
\label{determinants1}
\end{eqnarray}%
Therefore,$\vspace{-0.2in}$%
\begin{equation*}
\int\limits_{\mathbb{R}^{r^{2}+r}}\!\!\overline{f}_{p,\theta }\left(
x\right) \mathrm{d}\Phi \left( x;0,V_{p}\right)
\!\!=\!\!\prod_{j=1}^{r}\prod_{s=1}^{j}\sqrt{1\!-\!4\left( \!\theta
_{pj}v_{pj}\!\right) \left( \!\theta _{ps}v_{ps}\!\right) c_{p}}\vspace{%
-0.2in}
\end{equation*}%
and, uniformly in $p$ for $p$ sufficiently large$,\vspace{-0.2in}$%
\begin{equation}
\lim_{M,M^{\prime }\rightarrow \infty }\!\left\{
\!\int\limits_{B_{M,M^{\prime }}}\!\!\overline{f}_{p,\theta }\left( x\right) 
\mathrm{d}\Phi \left( x;0,V_{p}\right) \!-\!\prod_{j=1}^{r}\prod_{s=1}^{j}%
\sqrt{1\!-\!4\left( \!\theta _{pj}v_{pj}\!\right) \left( \!\theta
_{ps}v_{ps}\!\right) c_{p}}\right\} \!\!=\!\!0.\vspace{-0.2in}
\label{MMintegralFinal}
\end{equation}%
Equations (\ref{IMMnow}), (\ref{JMM_third}), and (\ref{MMintegralFinal})
describe the behavior of $I_{p}^{M,M^{\prime }}\left( \Theta _{p},\Lambda
_{p}\right) $ for large $p,$ $M,$ and $M^{\prime }$.

Let us now turn to the analysis of $I_{p}\left( \Theta _{p},\Lambda
_{p}\right) -I_{p}^{M,M^{\prime }}\left( \Theta _{p},\Lambda _{p}\right) .$
Let $B_{M}$ be the event $\left\{ \left\vert \gamma _{p1}^{(j,s)}\right\vert
\leq M\text{ for all }j,s\leq r\right\} $, and let$\vspace{-0.2in}$%
\begin{equation*}
I_{p}^{M}\left( \Theta _{p},\Lambda _{p}\right) =\mathbb{E}_{\Lambda
_{p}}\left( \mathbf{1}\left\{ B_{M}\right\} \exp \left\{
p\sum_{j=1}^{r}\theta _{pj}\frac{\tilde{g}^{(j)\prime }\Lambda _{p}\tilde{g}%
^{(j)}}{\tilde{g}^{(j)\prime }\tilde{g}^{(j)}}\right\} \right) .\vspace{%
-0.2in}
\end{equation*}%
As explained in Guionnet and Maida (2005, p.455), $\gamma _{p1}^{(j,s)},$ $%
j,s=1,...,r$ are independent of $\tilde{g}^{(j)\prime }\Lambda \tilde{g}%
^{(j)}/\tilde{g}^{(j)\prime }\tilde{g}^{(j)},$ $j=1,...,r.$ Therefore,$%
\vspace{-0.2in}$%
\begin{equation*}
I_{p}^{M}\left( \Theta _{p},\Lambda _{p}\right) =\mathbb{E}_{\Lambda
_{p}}\left( \mathbf{1}\left\{ B_{M}\right\} \right) I_{p}\left( \Theta
_{p},\Lambda _{p}\right) =\left( 1-\mathbb{E}_{\Lambda _{p}}\left( \mathbf{1}%
\left\{ B_{M}^{c}\right\} \right) \right) I_{p}\left( \Theta _{p},\Lambda
_{p}\right) .\vspace{-0.2in}
\end{equation*}

Denoting again by $\mathbb{P}$ the centered standard Gaussian measure on $%
\mathbb{R}$, we have$\vspace{-0.2in}$%
\begin{equation*}
\mathbb{E}_{\Lambda _{p}}\left( \mathbf{1}\left\{ \left\vert \gamma
_{p1}^{(j,s)}\right\vert \geq M\right\} \right) =\int \mathbf{1}\left\{
\left\vert \gamma _{p1}^{(j,s)}\right\vert \geq M\right\}
\prod_{j=1}^{r}\prod_{i=1}^{p}\mathrm{d}\mathbb{P}\left( g_{i}^{(j)}\right) .%
\vspace{-0.2in}
\end{equation*}%
For $j\neq s$ and $\tau \in \left( -\frac{1}{2}\sqrt{p},\frac{1}{2}\sqrt{p}%
\right) ,\vspace{-0.2in}$%
\begin{eqnarray*}
\int e^{\tau \gamma _{p1}^{(j,s)}}\prod_{j=1}^{r}\prod_{i=1}^{p}\mathrm{d}%
\mathbb{P}\left( g_{i}^{(j)}\right) &=&\frac{1}{\left( 2\pi \right) ^{p}}%
\int e^{\tau \frac{1}{\sqrt{p}}g^{(j)\prime }g^{(s)}}e^{-\frac{1}{2}\left(
g^{(j)\prime }g^{(j)}+g^{(s)\prime }g^{(s)}\right) }\prod_{i=1}^{p}\left( 
\mathrm{d}g_{i}^{(j)}\mathrm{d}g_{i}^{(s)}\right) \\
&=&\left( 1-\frac{\tau ^{2}}{p}\right) ^{-\frac{p}{2}}\leq e^{2\tau ^{2}}.%
\vspace{-0.2in}
\end{eqnarray*}%
Therefore, using Chebyshev's inequality, for $j\neq s$ and $\tau \in \left( -%
\frac{1}{2}\sqrt{p},\frac{1}{2}\sqrt{p}\right) ,\vspace{-0.2in}$%
\begin{equation*}
\int \mathbf{1}\left\{ \gamma _{p1}^{(j,s)}\geq M\right\}
\prod_{j=1}^{r}\prod_{i=1}^{p}\mathrm{d}\mathbb{P}\left( g_{i}^{(j)}\right)
\leq \frac{e^{2\tau ^{2}}}{e^{M\tau }}.\vspace{-0.2in}
\end{equation*}%
Setting $\tau =M/4$ (here we assume that $M<2\sqrt{p}$), we get$\vspace{%
-0.2in}$%
\begin{equation*}
\int \mathbf{1}\left\{ \gamma _{p1}^{(j,s)}\geq M\right\}
\prod_{j=1}^{r}\prod_{i=1}^{p}\mathrm{d}\mathbb{P}\left( g_{i}^{(j)}\right)
\leq e^{-M^{2}/8}.\vspace{-0.2in}
\end{equation*}%
Similarly, we show that the same inequality holds when $\gamma _{p1}^{(j,s)}$
is replaced by $-\gamma _{p}^{(j,s)},$ and thus$\vspace{-0.2in}$%
\begin{equation}
\int \mathbf{1}\left\{ \left\vert \gamma _{p1}^{(j,s)}\right\vert \geq
M\right\} \prod_{j=1}^{r}\prod_{i=1}^{p}\mathrm{d}\mathbb{P}\left(
g_{i}^{(j)}\right) \leq 2e^{-M^{2}/8}.  \label{subGaussGms}
\end{equation}%
For $j=s,$ the same line of arguments yields$\vspace{-0.2in}$%
\begin{equation}
\int \mathbf{1}\left\{ \left\vert \gamma _{p}^{(j,j)}\right\vert \geq
M\right\} \prod_{j=1}^{r}\prod_{i=1}^{p}\mathrm{d}\mathbb{P}\left(
g_{i}^{(j)}\right) \leq 2e^{-M^{2}/16}.\vspace{-0.2in}  \label{subGaussGmm}
\end{equation}

Inequalities (\ref{subGaussGms}) and (\ref{subGaussGmm}) imply that $\mathbb{%
E}_{\Lambda _{p}}\left( \mathbf{1}\left\{ B_{M}^{c}\right\} \right) \leq
2r^{2}e^{-M^{2}/16},$ and therefore, for sufficiently large $p,\vspace{-0.2in%
}$%
\begin{equation}
I_{p}\left( \Theta _{p},\Lambda _{p}\right) \geq I_{p}^{M}\left( \Theta
_{p},\Lambda _{p}\right) \geq \left( 1-2r^{2}e^{-M^{2}/16}\right)
I_{p}\left( \Theta _{p},\Lambda _{p}\right) .\vspace{-0.2in}
\label{IMversusI}
\end{equation}%
Note that$\vspace{-0.2in}$%
\begin{equation}
I_{p}^{M}\left( \!\Theta _{p},\Lambda _{p}\!\right) \!=e^{p\sum_{j=1}^{r} 
\left[ \theta _{pj}v_{pj}-\frac{1}{2p}\sum_{i=1}^{p}\ln \left( 1\!+\!2\theta
_{pj}v_{pj}\!-\!2\theta _{pj}\lambda _{i}\right) \right] }\left(
J_{p}^{M,M^{\prime }}+J_{p}^{M,M^{\prime },\infty }\right) ,\vspace{-0.2in}
\label{IpM}
\end{equation}%
where$\vspace{-0.2in}$%
\begin{equation*}
J_{p}^{M,M^{\prime },\infty }=\int \!\mathbf{1}\left\{ B_{M}\backslash
B_{M,M^{\prime }}\right\} f_{p,\theta }\left( \gamma _{p}\right)
\prod_{j=1}^{r}\prod_{i=1}^{p}\mathrm{d}\mathbb{P}_{pi}^{(j)}\left(
g_{i}^{(j)}\right) .\vspace{-0.2in}
\end{equation*}%
We will now derive an upper bound for $J_{p}^{M,M^{\prime },\infty }.$

From the definition of $f_{p,\theta }\left( \gamma _{p}\right) ,$ we see
that there exist positive constants $\beta _{1}$ and $\beta _{2},$ which may
depend on $r,\varepsilon $ and $\eta $, such that for any $\theta _{pj}$
satisfying $2\theta _{pj}\in \Omega _{\varepsilon \eta },$ $j\leq r$ and for
sufficiently large $p,$ when $B_{M}$ holds,$\vspace{-0.2in}$%
\begin{equation*}
f_{p,\theta }\left( \gamma _{p}\right) \leq \exp \left\{ \beta
_{1}M\sum_{s,k=1}^{r}\left\vert \gamma _{p2}^{(k,s)}\right\vert +\beta
_{2}M^{2}\right\} \text{.}\vspace{-0.2in}
\end{equation*}%
Let $B_{M,M^{\prime }}^{(k,s)}=B_{M}\cap \left\{ \left\vert \gamma
_{p2}^{(k,s)}\right\vert =\max_{j,m\leq r}\left\vert \gamma
_{p2}^{(j,m)}\right\vert >M^{\prime }\right\} $. Clearly, $B_{M}\backslash
B_{M,M^{\prime }}=\bigcup_{k,s=1}^{r}B_{M}^{(k,s)}$. Therefore,$\vspace{%
-0.2in}$%
\begin{eqnarray*}
J_{p}^{M,M^{\prime },\infty } &\leq
&\sum\limits_{k,s=1}^{r}\int_{B_{M,M^{\prime }}^{(k,s)}}e^{\beta
_{1}Mr^{2}\left\vert \gamma _{p2}^{(k,s)}\right\vert +\beta
_{2}M^{2}}\prod_{j=1}^{r}\prod_{i=1}^{p}\mathrm{d}\mathbb{P}%
_{pi}^{(j)}\left( g_{i}^{(j)}\right) \\
&\leq &\sum\limits_{j,m=1}^{r}\int_{\left\vert \gamma
_{p2}^{(k,s)}\right\vert \geq M^{\prime }}e^{\beta _{1}Mr^{2}\left\vert
\gamma _{p2}^{(k,s)}\right\vert +\beta
_{2}M^{2}}\prod_{j=1}^{r}\prod_{i=1}^{p}\mathrm{d}\mathbb{P}%
_{pi}^{(j)}\left( g_{i}^{(j)}\right) .\vspace{-0.2in}
\end{eqnarray*}

First assume $k\neq s.$ Denote $\lambda _{pi}\left( 1-2\theta _{pk}\lambda
_{pi}+2\theta _{pk}v_{pk}\right) ^{-1/2}\left( 1-2\theta _{ps}\lambda
_{pi}+2\theta _{ps}v_{ps}\right) ^{-1/2}$ as $\tilde{\lambda}_{pi}$ and $%
\left( 1-2\theta _{pj}\lambda _{pi}+2\theta _{pj}v_{pj}\right)
^{1/2}g_{i}^{(j)}$ as $\tilde{g}_{i}^{(j)}.$ Note that, under $\mathbb{P}%
_{pi}^{(j)},$ $\tilde{g}_{i}^{(j)}$ is a standard normal random variable.
Further, as long as~$2\theta _{pj}\in \Omega _{\varepsilon \eta }$ for $%
j\leq r,$ $\tilde{\lambda}_{pi}$ considered as a function of $\lambda _{i}$
is continuous on $\lambda _{i}\in \limfunc{supp}\hat{F}_{p}^{\lambda }$ for
sufficiently large $p,$ a.s.. Hence, the empirical distribution of $\tilde{%
\lambda}_{i}$ converges. Moreover, $\tilde{\lambda}_{\max }=\max_{i=1,...,p}(%
\tilde{\lambda}_{pi})$ and $\tilde{\lambda}_{\min }=\min_{i=1,...,p}(\tilde{%
\lambda}_{pi})$ a.s. converge to finite real numbers. Now, for $\tau $ such
that $\left\vert \tau \right\vert <\sqrt{p}/(2\tilde{\lambda}_{\max }),$ we
have$\vspace{-0.2in}$%
\begin{eqnarray*}
&&\int e^{\tau \gamma _{p2}^{(k,s)}}\prod_{j=1}^{r}\prod_{i=1}^{p}\mathrm{d}%
\mathbb{P}_{pi}^{(j)}\left( g_{i}^{(j)}\right) =\mathbb{E}e^{\tau \sqrt{p}%
\frac{1}{p}\sum_{i=1}^{p}\tilde{\lambda}_{pi}\tilde{g}_{i}^{\left( k\right) }%
\tilde{g}_{i}^{(s)}} \\
&=&\prod_{i=1}^{p}\mathbb{E}e^{\tau \frac{1}{\sqrt{p}}\tilde{\lambda}_{pi}%
\tilde{g}_{i}^{\left( k\right) }\tilde{g}_{i}^{(s)}}=\prod_{i=1}^{p}\left(
1-\tau ^{2}\frac{\tilde{\lambda}_{pi}^{2}}{p}\right) ^{-1/2}\leq e^{2\tilde{%
\lambda}_{\max }^{2}\tau ^{2}}\vspace{-0.2in}
\end{eqnarray*}%
for sufficiently large $p,$ a.s.. Using this inequality, we get, for
sufficiently large $p$ and any positive $t$ such that $\beta _{1}r^{2}M+t<%
\sqrt{p}/(2\tilde{\lambda}_{\max }),\vspace{-0.2in}$%
\begin{eqnarray*}
&&\int_{\gamma _{p2}^{(k,s)}\geq M^{\prime }}\!e^{\beta _{1}r^{2}M\gamma
_{p2}^{(k,s)}}\!\prod_{j=1}^{r}\prod_{i=1}^{p}\mathrm{d}\mathbb{P}%
_{pi}^{(j)}\left( \!g_{i}^{(j)}\!\right) \!\leq \!\int \!e^{\beta
_{1}r^{2}M\gamma _{p2}^{(k,s)}\!+\!t\left( \gamma
_{p2}^{(k,s)}\!-\!M^{\prime }\right) }\!\prod_{j=1}^{r}\prod_{i=1}^{p}%
\mathrm{d}\mathbb{P}_{pi}^{(j)}\left( \!g_{i}^{(j)}\!\right) \\
&=&e^{-tM^{\prime }}\int \!e^{\left( \beta _{1}r^{2}M+t\right) \gamma
_{p2}^{(k,s)}}\!\prod_{j=1}^{r}\prod_{i=1}^{p}\mathrm{d}\mathbb{P}%
_{pi}^{(j)}\left( g_{i}^{(j)}\right) \!\leq \!e^{-tM^{\prime }}e^{2\tilde{%
\lambda}_{\max }^{2}\left( \beta _{1}r^{2}M+t\right) ^{2}}.\vspace{-0.2in}
\end{eqnarray*}%
Setting $t=\frac{M^{\prime }}{4\tilde{\lambda}_{\max }^{2}}-\beta _{1}r^{2}M$
(here we assume that $M$ and $M^{\prime }$ are such that $t$ satisfies the
above requirements), we get$\vspace{-0.2in}$%
\begin{equation*}
\int_{\gamma _{p2}^{(k,s)}\geq M^{\prime }}e^{\beta _{1}r^{2}M\gamma
_{p2}^{(k,s)}}\prod_{j=1}^{r}\prod_{i=1}^{p}\mathrm{d}\mathbb{P}%
_{pi}^{(j)}\left( \!g_{i}^{(j)}\!\right) \leq e^{-\frac{\left( M^{\prime
}\right) ^{2}}{8\tilde{\lambda}_{\max }^{2}}+\beta _{1}r^{2}MM^{\prime }}.%
\vspace{-0.2in}
\end{equation*}%
Replacing $\gamma _{p2}^{(k,s)}$ by $-\gamma _{p2}^{(k,s)}$ in the above
derivations and combining the result with the above inequality, we get$%
\vspace{-0.2in}$%
\begin{equation*}
\int_{\left\vert \gamma _{p2}^{(k,s)}\right\vert \geq M^{\prime }}e^{\beta
_{1}r^{2}M\left\vert \gamma _{p2}^{(k,s)}\right\vert
}\prod_{j=1}^{r}\prod_{i=1}^{p}\mathrm{d}\mathbb{P}_{pi}^{(j)}\left(
\!g_{i}^{(j)}\!\right) \leq 2e^{-\frac{\left( M^{\prime }\right) ^{2}}{8%
\tilde{\lambda}_{\max }^{2}}+\beta _{1}r^{2}MM^{\prime }}.\vspace{-0.2in}
\end{equation*}%
When $k=s,$ following a similar line of arguments, we obtain$\vspace{-0.2in}$%
\begin{equation*}
\int_{\left\vert \gamma _{p2}^{(k,k)}\right\vert \geq M^{\prime }}e^{\beta
_{1}r^{2}M\left\vert \gamma _{p2}^{(k,k)}\right\vert
}\prod_{j=1}^{r}\prod_{i=1}^{p}\mathrm{d}\mathbb{P}_{pi}^{(j)}\left(
\!g_{i}^{(j)}\!\right) \leq 4e^{-\frac{\left( M^{\prime }\right) ^{2}}{16%
\tilde{\lambda}_{\max }^{2}}+\beta _{1}r^{2}MM^{\prime }}.\vspace{-0.2in}
\end{equation*}%
and thus, for sufficiently large $p,\vspace{-0.2in}$%
\begin{equation}
J_{p}^{M,M^{\prime },\infty }\leq 4r^{2}e^{-\frac{\left( M^{\prime }\right)
^{2}}{16\tilde{\lambda}_{\max }^{2}}+\beta _{1}r^{2}MM^{\prime }}.\vspace{%
-0.2in}  \label{Jkk}
\end{equation}

Finally, combining (\ref{IMversusI}), (\ref{IpM}), and (\ref{Jkk}), we
obtain for$\vspace{-0.2in}$%
\begin{equation}
J_{p}=I_{p}\left( \Theta _{p},\Lambda _{p}\right) e^{-p\sum_{j=1}^{r}\left[
\theta _{pj}v_{pj}-\frac{1}{2p}\sum_{i=1}^{p}\ln \left( 1\!+\!2\theta
_{pj}v_{pj}\!-\!2\theta _{pj}\lambda _{i}\right) \right] }\vspace{-0.2in}
\label{Jp}
\end{equation}
the following upper and lower bounds:%
\begin{equation}
J_{p}^{M,M^{\prime }}\leq J_{p}\leq \left( 1-2r^{2}e^{-\frac{M^{2}}{16}%
}\right) ^{-1}\left( J_{p}^{M,M^{\prime }}+4r^{2}e^{-\frac{\left( M^{\prime
}\right) ^{2}}{16\tilde{\lambda}_{\max }^{2}}+\beta _{1}r^{2}MM^{\prime
}}\right) .\vspace{-0.2in}  \label{lastIneq}
\end{equation}%
Let $\tau >0$ be an arbitrarily small number. Equations (\ref{JMM_third})
and (\ref{MMintegralFinal}) imply that there exist $\bar{M}$ and $\bar{M}%
^{\prime }$ such that, for any $M>\bar{M}$ and $M^{\prime }>\bar{M}^{\prime
},\vspace{-0.2in}$%
\begin{equation*}
\left\vert J_{p}^{M,M^{\prime }}-\prod_{j=1}^{r}\prod_{s=1}^{j}\sqrt{%
1\!-\!4\left( \!\theta _{pj}v_{pj}\!\right) \left( \!\theta
_{ps}v_{ps}\!\right) c_{p}}\right\vert <\frac{\tau }{4}\vspace{-0.2in}
\end{equation*}%
for all sufficiently large $p.$ Let us choose $M>\bar{M}$ and $M^{\prime }>%
\bar{M}^{\prime }$ so that$\vspace{-0.2in}$%
\begin{eqnarray*}
\left( 1-2r^{2}e^{-\frac{M^{2}}{16}}\right) ^{-1} &<&2, \\
\left( 1-2r^{2}e^{-\frac{M^{2}}{16}}\right) ^{-1}4r^{2}e^{-\frac{\left(
M^{\prime }\right) ^{2}}{16\tilde{\lambda}_{\max }^{2}}+\beta
_{1}r^{2}MM^{\prime }} &<&\frac{\tau }{4},\vspace{-0.2in}
\end{eqnarray*}%
and$\vspace{-0.2in}$%
\begin{equation*}
\left[ \left( 1-2r^{2}e^{-\frac{M^{2}}{16}}\right) ^{-1}-1\right]
\sup_{\left\{ 2\theta _{pj}\in \Omega _{\varepsilon \eta },j\leq r\right\}
}\prod_{j=1}^{r}\prod_{s=1}^{j}\sqrt{1\!-\!4\left( \!\theta
_{pj}v_{pj}\!\right) \left( \!\theta _{ps}v_{ps}\!\right) c_{p}}<\frac{\tau 
}{4}\vspace{-0.2in}
\end{equation*}%
for all sufficiently large $p,$ a.s.. Then, (\ref{lastIneq}) implies that$%
\vspace{-0.2in}$%
\begin{equation}
\left\vert J_{p}-\prod_{j=1}^{r}\prod_{s=1}^{j}\sqrt{1\!-\!4\left( \!\theta
_{pj}v_{pj}\!\right) \left( \!\theta _{ps}v_{ps}\!\right) c_{p}}\right\vert
<\tau \vspace{-0.2in}  \label{Jpbound}
\end{equation}%
for all sufficiently large $p$, a.s.. Since $\tau $ can be chosen
arbitrarily, we have, from (\ref{Jp}) and (\ref{Jpbound}),$\vspace{-0.2in}$%
\begin{eqnarray*}
I_{p}\left( \Theta _{p},\Lambda _{p}\right) &=&e^{p\sum_{j=1}^{r}\left[
\theta _{pj}v_{pj}-\frac{1}{2p}\sum_{i=1}^{p}\ln \left( 1\!+\!2\theta
_{pj}v_{pj}\!-\!2\theta _{pj}\lambda _{pi}\right) \right] } \\
&&\times \left( \prod_{j=1}^{r}\prod_{s=1}^{j}\sqrt{1\!-\!4\left( \!\theta
_{pj}v_{pj}\!\right) \left( \!\theta _{ps}v_{ps}\!\right) c_{p}}+o(1)\right)
,\vspace{-0.2in}
\end{eqnarray*}%
where $o(1)\rightarrow 0$ as $p\rightarrow \infty $ uniformly in $\left\{
2\theta _{pj}\in \Omega _{\varepsilon \eta },j\leq r\right\} ,$ a.s..$%
\square $

\subsection{Proof of Theorem \protect\ref{Theorem1}}

Setting $\theta _{pj}=\frac{1}{2c_{p}}\frac{h_{j}}{1+h_{j}}$, we have $%
v_{pj}=1+h_{j}$, $\theta _{pj}v_{pj}=\frac{h_{j}}{2c_{p}}$, and\vspace{-0.2in%
}%
\begin{equation*}
\ln \left( 1\!+\!2\theta _{pj}v_{pj}\!-\!2\theta _{pj}\lambda _{pi}\right)
=\ln \left( \frac{1}{c_{p}}\frac{h_{j}}{1+h_{j}}\right) +\ln \left(
z_{0j}-\lambda _{pi}\right) .\vspace{-0.2in}
\end{equation*}
Further, by Lemma 11 and formula (3.3) of OMH, $\int \ln \left(
z_{j0}-\lambda \right) \mathrm{d}F_{p}^{MP}\left( \lambda \right) =\frac{%
h_{j}}{c_{p}}-\frac{1}{c_{p}}\ln \left( 1+h_{j}\right) +\ln \frac{\left(
1+h_{j}\right) c_{p}}{h_{j}}$ for sufficiently large $p,$ a.s.. With these
auxiliary results, formula (\ref{equivalence 1}) is a straightforward
consequence of (\ref{LR1}) and Proposition 2.

Turning to the proof of (\ref{equivalence 2}), consider the integrals$%
\vspace{-0.2in}$

\begin{equation*}
\mathcal{I}\left( k_{1},k_{2}\right) \!=\!\int_{k_{1}}^{k_{2}}\!\!\!x^{\frac{%
n_{p}p}{2}-1}e^{-\frac{n_{p}}{2}x}\!\!\int\limits_{\mathcal{O}\left(
p\right) }\!\!\!e^{p\frac{x}{S_{p}}\limfunc{tr}\left( D_{p}Q^{\prime
}\Lambda _{p}Q\right) }\left( \!\mathrm{d}Q\!\right) \!\mathrm{d}x,\text{ }%
k_{1}<k_{2}\in \mathbb{R}.\vspace{-0.2in}
\end{equation*}%
In what follows, we omit the subscript $p$ in $n_{p}$ to simplify notation.
Note that~$\mathcal{I}\left( 0,\infty \right) $ is the integral appearing in
expression (\ref{LR2}) for $L_{p}\!\left( h;\mu _{p}\right) $. Let us now
prove that, for some constant $\alpha >0,\vspace{-0.2in}$%
\begin{equation}
\mathcal{I}\left( 0,\infty \right) \!=\!\mathcal{I}\left( p\!-\!\alpha \sqrt{%
p},p\!+\!\alpha \sqrt{p}\right) \left( 1\!+\!o\left( 1\right) \right) ,\text{
a.s.}\vspace{-0.2in}  \label{IP}
\end{equation}%
where $o(1)$ is uniform in $h\in \left[ 0,\sqrt{c}-\delta \right] ^{r}.$

Since, by Corollary 1, $S_{p}/p\rightarrow 1$ a.s., the set $H_{\delta }$ is
bounded from below, and $\lambda _{p1}\rightarrow \left( 1+\sqrt{c}\right)
^{2}$ a.s., there exists a constant $A_{1}>0$ that depends only on $\delta $
and $r,$ such that $\inf_{\left[ 0,\sqrt{c}-\delta \right] ^{r}}px\limfunc{tr%
}\left( D_{p}Q^{\prime }\Lambda _{p}Q\right) /S_{p}\geq -A_{1}x/2$ for all $%
x\geq 0$ and all sufficiently large $p,$ a.s.. Therefore, for all $h\in %
\left[ 0,\sqrt{c}-\delta \right] ^{r},\vspace{-0.2in}$%
\begin{equation*}
2\mathcal{I}\left( 0,\infty \right) \geq \int_{0}^{\infty }\!\!x^{\frac{np}{2%
}-1}e^{-\frac{n+A_{1}}{2}x}\mathrm{d}x=\left( \frac{n+A_{1}}{2}\right) ^{-%
\frac{np}{2}}\Gamma \left( \frac{np}{2}\right) ,\vspace{-0.2in}
\end{equation*}%
and, using Stirling's approximation, we get$\vspace{-0.2in}$%
\begin{eqnarray}
\mathcal{I}\left( 0,\infty \right) &\geq &\left( \frac{n+A_{1}}{2}\right) ^{-%
\frac{np}{2}}\left( \frac{np}{2}\right) ^{\frac{np}{2}}e^{-\frac{np}{2}%
}\left( \frac{4\pi }{np}\right) ^{1/2}\left( 1+o\left( 1\right) \right) 
\notag \\
&=&p^{\frac{np}{2}}e^{-\left( \frac{n}{2}+\frac{A_{1}}{2}-\frac{1}{4}\frac{%
A_{1}^{2}}{n}\right) p}\left( \frac{4\pi }{np}\right) ^{1/2}\left( 1+o\left(
1\right) \right) ,\text{ a.s.}\vspace{-0.2in}  \label{I1bound}
\end{eqnarray}

Next, there exists a constant $A_{2}>0$ such that, for all $x\geq 0$ and all
sufficiently large $p,$ $\sup_{h\in \left[ 0,\sqrt{c}-\delta \right] ^{r}}px%
\limfunc{tr}\left( D_{p}Q^{\prime }\Lambda _{p}Q\right) /S_{p}\leq A_{2}x/2$%
, a.s.. Therefore, a.s., for all sufficiently large $p$,$\vspace{-0.2in}$%
\begin{equation*}
\mathcal{I}\left( p\!+\!\alpha \sqrt{p},\infty \right) \!\leq
\!\int_{p+\alpha \sqrt{p}}^{\infty }\!\!x^{\frac{np}{2}-1}e^{-\frac{n-A_{2}}{%
2}x}\mathrm{d}x=\!\left( \frac{n-A_{2}}{2}\right) ^{-\frac{np}{2}}\Gamma
\left( \frac{np}{2},y\right) ,\vspace{-0.2in}
\end{equation*}%
where $\Gamma \left( \frac{np}{2},y\right) $ is the complementary incomplete
Gamma function (see Olver 1997, p.45) with $y=\left( p+\alpha \sqrt{p}%
\right) \left( \frac{n-A_{2}}{2}\right) .$ Hence, for sufficiently large $p,$
$y>np/2+n\alpha \sqrt{p}/4,$ and we can continue$\vspace{-0.2in}$%
\begin{equation*}
\mathcal{I}\left( p\!+\!\alpha \sqrt{p},\infty \right) \!<\left( \frac{%
n-A_{2}}{2}\right) ^{-\frac{np}{2}}\Gamma \left( \frac{np}{2},\frac{np}{2}+%
\frac{\alpha n\sqrt{p}}{4}\right) ,\text{ a.s.}\vspace{-0.2in}
\end{equation*}%
Now, $\Gamma \left( \beta ,\gamma \right) \leq e^{-\gamma }\gamma ^{\beta
}/(\gamma -\beta +1)$ whenever $\beta >1$ and $\gamma >\beta -1$ (Olver
1997, p.70). Therefore, we have, for sufficiently large $p$,$\vspace{-0.2in}$%
\begin{eqnarray*}
\mathcal{I}\left( p\!+\!\alpha \sqrt{p},\infty \right) \! &<&\!\left( 1-%
\frac{A_{2}}{n}\right) ^{-\frac{np}{2}}\frac{e^{-\frac{np}{2}\!-\!\frac{%
\alpha n\sqrt{p}}{4}}p^{\frac{np}{2}}\left( 1+\frac{\alpha }{2\sqrt{p}}%
\right) ^{\frac{np}{2}}}{\alpha n\sqrt{p}/4+1} \\
&=&p^{\frac{np}{2}}e^{\frac{A_{2}p}{2}+\frac{A_{2}^{2}p}{4n}}\frac{e^{-\frac{%
np}{2}-\frac{\alpha ^{2}n}{16}+\frac{\alpha ^{3}n}{48\sqrt{p}}-\frac{\alpha
^{4}n}{128p}}}{\alpha n\sqrt{p}/4+1}\left( 1+o(1)\right) \\
&<&p^{\frac{np}{2}}e^{-\frac{np}{2}}\frac{e^{p\left( A_{2}-\frac{\alpha ^{2}n%
}{32p}\right) }}{\alpha n\sqrt{p}/4+1}\left( 1+o(1)\right) ,\text{ a.s.}%
\vspace{-0.2in}
\end{eqnarray*}%
Comparing this to (\ref{I1bound}), we see that $\alpha $ can be chosen so
that$\vspace{-0.2in}$%
\begin{equation}
\mathcal{I}\left( p\!+\!\alpha \sqrt{p},\infty \right) =o(1)\mathcal{I}%
\left( 0,\infty \right) ,\text{ a.s.}\vspace{-0.2in}  \label{I2bound}
\end{equation}

Further, for sufficiently large $p$,$\vspace{-0.2in}$%
\begin{eqnarray*}
\mathcal{I}\left( 0,p\!-\!\alpha \sqrt{p}\right) &\leq &\int_{0}^{p-\alpha 
\sqrt{p}}\!\!x^{\frac{np}{2}-1}e^{-\frac{n-A_{2}}{2}x}\mathrm{d}x \\
&=&\left( \frac{n-A_{2}}{2}\right) ^{-\frac{np}{2}}\int_{0}^{y}\!\!t^{\frac{%
np}{2}-1}e^{-t}\mathrm{d}t,\text{ a.s.,}\vspace{-0.2in}
\end{eqnarray*}%
where $y=\left( p-\alpha \sqrt{p}\right) \frac{n-A_{2}}{2}<\frac{np}{2}-%
\frac{\alpha n\sqrt{p}}{4}$. Therefore, for any positive $z<\frac{np}{2}$
and sufficiently large $p,\vspace{-0.2in}$%
\begin{eqnarray*}
\mathcal{I}\left( 0,p\!-\!\alpha \sqrt{p}\right) &\leq &\left( \frac{%
n\!-\!A_{2}}{2}\right) ^{-\frac{np}{2}}\int_{0}^{\frac{np}{2}-\frac{\alpha n%
\sqrt{p}}{4}}\!\!t^{\frac{np}{2}-1}e^{-t}\mathrm{d}t \\
&<&\left( \frac{n\!-\!A_{2}}{2}\right) ^{-\frac{np}{2}}\left( \frac{np}{2}%
\!-\!\frac{\alpha n\sqrt{p}}{4}\right) ^{z}\Gamma \left( \frac{np}{2}%
\!-\!z\right) \text{.}\vspace{-0.2in}
\end{eqnarray*}%
Setting $z=\alpha n\sqrt{p}/4$ and using Stirling's approximation, we have,
a.s.,$\vspace{-0.2in}$%
\begin{equation*}
\left( \frac{np}{2}\!-\!\frac{\alpha n\sqrt{p}}{4}\right) ^{z}\Gamma \left( 
\frac{np}{2}\!-\!z\right) =\left( \frac{np}{2}\!-\!\frac{\alpha n\sqrt{p}}{4}%
\right) ^{\frac{np}{2}-\frac{1}{2}}e^{-\frac{np}{2}+\frac{\alpha n\sqrt{p}}{4%
}}\sqrt{2\pi }\left( 1\!+\!o\left( 1\right) \right) \vspace{-0.2in}
\end{equation*}%
so that$\vspace{-0.2in}$%
\begin{eqnarray*}
\mathcal{I}\left( 0,p\!-\!\alpha \sqrt{p}\right) &<&\left( \frac{n\!-\!A_{2}%
}{2}\right) ^{-\frac{np}{2}}\left( \frac{np}{2}\!-\!\frac{\alpha n\sqrt{p}}{4%
}\right) ^{\frac{np}{2}-\frac{1}{2}}e^{-\frac{np}{2}\!+\!\frac{\alpha n\sqrt{%
p}}{4}}\sqrt{2\pi }\left( 1\!+\!o\left( 1\right) \right) \\
&<&p^{\frac{np}{2}}e^{-\frac{np}{2}}e^{p\left( \frac{A_{2}}{2}+\frac{%
A_{2}^{2}}{4n}-\frac{\alpha ^{2}n}{16p}\right) }\left( 1\!+\!o\left(
1\right) \right) \text{, a.s..}\vspace{-0.2in}
\end{eqnarray*}%
Comparing this to (\ref{I1bound}), we see that $\alpha $ can be chosen so
that\vspace{-0.2in}%
\begin{equation}
\mathcal{I}\left( 0,p\!-\!\alpha \sqrt{p}\right) =o(1)\mathcal{I}\left(
0,\infty \right) ,\vspace{-0.2in}  \label{I3bound}
\end{equation}%
a.s.. Combining (\ref{I2bound}) and (\ref{I3bound}), we get (\ref{IP}).

Now, letting $\tilde{\theta}_{pj}=\frac{x}{S_{p}}\theta _{pj}=\frac{x}{S_{p}}%
\frac{1}{2c_{p}}\frac{h_{j}}{1+h_{j}},$ note that there exist $\varepsilon
>0 $ and $\eta >0$ such that $\left\{ 2\tilde{\theta}_{pj}:h_{j}\in \left[ 0,%
\sqrt{c}-\delta \right] \text{ and }x\in \left[ p\!-\!\alpha \sqrt{p}%
,p\!+\!\alpha \sqrt{p}\right] \right\} \subseteq \Theta _{\varepsilon \eta }$
for all sufficiently large $p,$ a.s.. Hence, by (\ref{IP}), and Proposition
2, a.s.,$\vspace{-0.2in}$%
\begin{eqnarray}
\mathcal{I}\left( 0,\infty \right) \!\! &=&\!\!\!\int_{p\!-\!\alpha \sqrt{p}%
}^{p\!+\!\alpha \sqrt{p}}\!\!\!x^{\frac{np}{2}-1}e^{-\frac{n}{2}%
x}e^{p\sum_{j=1}^{r}\left[ \tilde{\theta}_{pj}\tilde{v}_{pj}-\frac{1}{2p}%
\sum_{i=1}^{p}\ln \left( 1\!+\!2\tilde{\theta}_{pj}\tilde{v}_{pj}\!-\!2%
\tilde{\theta}_{pj}\lambda _{pi}\right) \right] }  \label{I0inf} \\
&&\times \left( \prod_{j=1}^{r}\prod_{s=1}^{j}\sqrt{1\!-\!4\left( \!\tilde{%
\theta}_{pj}\tilde{v}_{pj}\!\right) \left( \!\tilde{\theta}_{ps}\tilde{v}%
_{ps}\!\right) c_{p}}+o(1)\right) \!\mathrm{d}x\text{,}\vspace{-0.2in} 
\notag
\end{eqnarray}%
where $o(1)$ is uniform in $h\in \left[ 0,\sqrt{c}-\delta \right] ^{r}$ and $%
x\in \left[ p\!-\!\alpha \sqrt{p},p\!+\!\alpha \sqrt{p}\right] .$

Expanding $\tilde{\theta}_{pj}\tilde{v}_{pj}-\frac{1}{2p}\sum_{i=1}^{p}\ln
\left( 1\!+\!2\tilde{\theta}_{pj}\tilde{v}_{pj}\!-\!2\tilde{\theta}%
_{pj}\lambda _{pi}\right) $ and $\left( \!\tilde{\theta}_{pj}\tilde{v}%
_{pj}\!\right) \left( \!\tilde{\theta}_{ps}\tilde{v}_{ps}\!\right) $ into
power series of $\frac{x}{p}-1,$ we get$\vspace{-0.2in}$%
\begin{eqnarray*}
\mathcal{I}\left( 0,\infty \right) \!\! &=&\!\!\int_{p\!-\!\alpha \sqrt{p}%
}^{p\!+\!\alpha \sqrt{p}}\!\!\!x^{\frac{np}{2}-1}e^{-\frac{n}{2}x}e^{p\left(
B_{0}+B_{1}\left( \frac{x}{p}-1\right) +B_{2}\left( \frac{x}{p}-1\right)
^{2}\right) } \\
&&\times \left( \prod_{j=1}^{r}\prod_{s=1}^{j}\sqrt{1\!-\!4\left( \!\theta
_{pj}v_{pj}\!\right) \left( \!\theta _{ps}v_{ps}\!\right) c_{p}}+o(1)\right)
\!\mathrm{d}x\text{,}\vspace{-0.2in}
\end{eqnarray*}%
where $B_{0},B_{1}$ and $B_{2}$ are $O(1)$ uniformly in $h\in \left[ 0,\sqrt{%
c}-\delta \right] ^{r}.$ Further, consider the integral$\vspace{-0.2in}$%
\begin{equation*}
I^{(0)}=\int_{p\!-\!\alpha \sqrt{p}}^{p\!+\!\alpha \sqrt{p}}\!\!\!x^{\frac{np%
}{2}-1}e^{-\frac{n}{2}x}\!\!e^{p\left( B_{1}\frac{x}{p}\!+\!B_{2}\left( 
\frac{x}{p}-1\right) ^{2}\right) }\mathrm{d}x.\vspace{-0.2in}
\end{equation*}%
Splitting the domain of integration into segments $\left[ p\!-\!\alpha \sqrt{%
p},p\!-\!\alpha p^{\gamma }\right] ,\left[ p\!-\!\alpha p^{\gamma
},p\!+\!\alpha p^{\gamma }\right] $ and $\left[ p\!+\!\alpha p^{\gamma
},p\!+\!\alpha \sqrt{p}\right] ,$ where $0<\gamma <1/2,$ and denoting the
corresponding integrals by $I^{(1)},I^{(2)}$ and $I^{(3)},$ respectively, we
have$\vspace{-0.2in}$%
\begin{eqnarray*}
I^{(1)} &<&e^{\alpha ^{2}}\int_{p\!-\!\alpha \sqrt{p}}^{p\!-\!\alpha
p^{\gamma }}\!\!\!x^{\frac{np}{2}-1}e^{-\frac{n}{2}x}e^{B_{1}x\!}\mathrm{d}%
x<e^{\alpha ^{2}}p^{\frac{np}{2}}\left( 1-\frac{2B_{1}}{n}\right) ^{\frac{np%
}{2}}\int_{0}^{1\!-\!\frac{\alpha }{2}p^{\gamma -1}}\!\!\!y^{\frac{np}{2}%
-1}e^{-\frac{np}{2}y}\mathrm{d}y, \\
I^{(2)} &>&\int_{p\!-\!\alpha p^{\gamma }}^{p\!+\!\alpha p^{\gamma
}}\!\!\!x^{\frac{np}{2}-1}e^{-\frac{n}{2}x}e^{B_{1}x\!}\mathrm{d}x>p^{\frac{%
np}{2}}\left( 1-\frac{2B_{1}}{n}\right) ^{\frac{np}{2}}\int_{1\!-\!\frac{%
\alpha }{2}p^{\gamma -1}}^{1\!+\!\frac{\alpha }{2}p^{\gamma -1}}\!\!\!y^{%
\frac{np}{2}-1}e^{-\frac{np}{2}y}\mathrm{d}y,\text{ and} \\
I^{(3)} &<&e^{\alpha ^{2}}\int_{p\!+\!\alpha p^{\gamma }}^{p\!+\!\alpha 
\sqrt{p}}\!\!\!x^{\frac{np}{2}-1}e^{-\frac{n}{2}x}e^{B_{1}x\!}\mathrm{d}%
x<e^{\alpha ^{2}}p^{\frac{np}{2}}\left( 1-\frac{2B_{1}}{n}\right) ^{\frac{np%
}{2}}\int_{1\!+\!\frac{\alpha }{2}p^{\gamma -1}}^{\infty }\!\!\!y^{\frac{np}{%
2}-1}e^{-\frac{np}{2}y}\mathrm{d}y.\vspace{-0.2in}
\end{eqnarray*}%
Using the Laplace approximation, we have$\vspace{-0.2in}$%
\begin{eqnarray*}
\int_{0}^{1\!-\!\frac{\alpha }{2}p^{\gamma -1}}\!\!\!y^{\frac{np}{2}-1}e^{-%
\frac{np}{2}y}\mathrm{d}y &=&o(1)\int_{1\!-\!\frac{\alpha }{2}p^{\gamma
-1}}^{1\!+\!\frac{\alpha }{2}p^{\gamma -1}}\!\!\!y^{\frac{np}{2}-1}e^{-\frac{%
np}{2}y}\mathrm{d}y\text{, and} \\
\int_{1\!+\!\frac{\alpha }{2}p^{\gamma -1}}^{\infty }\!\!\!y^{\frac{np}{2}%
-1}e^{-\frac{np}{2}y}\mathrm{d}y &=&o(1)\int_{1\!-\!\frac{\alpha }{2}%
p^{\gamma -1}}^{1\!+\!\frac{\alpha }{2}p^{\gamma -1}}\!\!\!y^{\frac{np}{2}%
-1}e^{-\frac{np}{2}y}\mathrm{d}y,\vspace{-0.2in}
\end{eqnarray*}%
so that $I^{(2)}$ dominates $I^{(1)}$ and $I^{(3)}$ and$\vspace{-0.2in}$%
\begin{eqnarray*}
I^{(0)} &=&\left( 1+o(1)\right) \int_{p\!-\!\alpha p^{\gamma
}}^{p\!+\!\alpha p^{\gamma }}\!\!\!x^{\frac{np}{2}-1}e^{-\frac{n}{2}%
x}e^{p\left( B_{1}\frac{x}{p}\!+\!B_{2}\left( \frac{x}{p}-1\right)
^{2}\right) }\mathrm{d}x \\
&=&\left( 1+o(1)\right) \int_{p\!-\!\alpha p^{\gamma }}^{p\!+\!\alpha
p^{\gamma }}\!\!\!x^{\frac{np}{2}-1}e^{-\frac{n}{2}x}e^{B_{1}x\!}\mathrm{d}x
\\
&=&\left( 1+o(1)\right) \int_{p\!-\!\alpha \sqrt{p}}^{p\!+\!\alpha \sqrt{p}%
}\!\!\!x^{\frac{np}{2}-1}e^{-\frac{n}{2}x}e^{B_{1}x\!}\mathrm{d}x.\vspace{%
-0.2in}
\end{eqnarray*}%
This implies that$\vspace{-0.2in}$%
\begin{eqnarray*}
\mathcal{I}\left( 0,\infty \right) \!\! &=&\!\!\int_{p\!-\!\alpha \sqrt{p}%
}^{p\!+\!\alpha \sqrt{p}}\!\!\!x^{\frac{np}{2}-1}e^{-\frac{n}{2}x}e^{p\left(
B_{0}+B_{1}\left( \frac{x}{p}-1\right) \right) } \\
&&\times \left( \prod_{j=1}^{r}\prod_{s=1}^{j}\sqrt{1\!-\!4\left( \!\theta
_{pj}v_{pj}\!\right) \left( \!\theta _{ps}v_{ps}\!\right) c_{p}}+o(1)\right)
\!\mathrm{d}x\text{,}\vspace{-0.2in}
\end{eqnarray*}%
and hence, only constant and linear terms in the expansion of\linebreak\ $%
\tilde{\theta}_{pj}\tilde{v}_{pj}-\frac{1}{2p}\sum_{i=1}^{p}\ln \left(
1\!+\!2\tilde{\theta}_{pj}\tilde{v}_{pj}\!-\!2\tilde{\theta}_{pj}\lambda
_{pi}\right) $ into power series of $\frac{x}{p}-1$ matter for the
evaluation of $\mathcal{I}\left( 0,\infty \right) .$ Let us find these terms.

By Corollary 1, $\frac{x}{S_{p}}-1=\frac{x}{p}-\frac{S_{p}}{p}+o(p^{-1})$
a.s.. Using this fact, after some algebra, we get$\vspace{-0.2in}$%
\begin{equation*}
\tilde{\theta}_{pj}\widetilde{v}_{pj}\!=\!\theta _{pj}v_{pj}+\theta
_{pj}v_{pj}^{2}\left( \frac{x}{p}-\frac{S_{p}}{p}\right) +O\left( \left( 
\frac{x}{p}-1\right) ^{2}\right) ,\vspace{-0.2in}
\end{equation*}%
\begin{equation*}
\ln \left( 2\tilde{\theta}_{pj}\right) \!=\!\ln \left( 2\theta _{pj}\right)
\!+\!\left( \frac{x}{p}\!-\!\frac{S_{p}}{p}\right) +O\left( \left( \frac{x}{p%
}\!-\!1\right) ^{2}\right) ,\vspace{-0.2in}
\end{equation*}%
and$\vspace{-0.2in}$%
\begin{eqnarray*}
\sum_{i=1}^{p}\ln \left( K_{p}^{MP}\left( 2\tilde{\theta}_{pj}\right)
\!-\!\lambda _{pi}\right) \! &=&\!\sum_{i=1}^{p}\ln \left( K_{p}^{MP}\left(
2\theta _{pj}\right) \!-\!\lambda _{pi}\right) -p\left( 1\!-\!4c_{p}\theta
_{pj}^{2}v_{pj}^{2}\right) \left( \frac{x}{p}\!-\!\frac{S_{p}}{p}\right) \\
&&+O\left( \left( \frac{x}{p}-1\right) ^{2}\right) .\vspace{-0.2in}
\end{eqnarray*}%
It follows that$\vspace{-0.2in}$%
\begin{eqnarray}
\mathcal{I}\left( 0,\infty \right) \!\! &=&\!\!\int_{p\!-\!\alpha \sqrt{p}%
}^{p\!+\!\alpha \sqrt{p}}\!\!\!x^{\frac{np}{2}-1}e^{-\frac{n}{2}%
x}e^{p\sum_{j=1}^{r}\left[ \theta _{pj}v_{pj}-\frac{1}{2p}\sum_{i=1}^{p}\ln
\left( 1\!+\!2\theta _{pj}v_{pj}\!-\!2\theta _{pj}\lambda _{pi}\right) %
\right] }  \label{last eq} \\
&&\times e^{\sum_{j=1}^{r}\theta _{pj}v_{pj}\left( x-S_{p}\right) }\left(
\prod_{j=1}^{r}\prod_{s=1}^{j}\sqrt{1\!-\!4\left( \!\theta
_{pj}v_{pj}\!\right) \left( \!\theta _{ps}v_{ps}\!\right) c_{p}}+o(1)\right)
\!\mathrm{d}x  \notag \\
&=&\left( 1+o(1)\right) \prod_{j=1}^{r}\left( 1+h_{j}\right) ^{\frac{n_{p}}{2%
}}L_{p}\!\left( h;\lambda _{p}\right) \int_{p\!-\!\alpha \sqrt{p}%
}^{p\!+\!\alpha \sqrt{p}}\!\!\!x^{\frac{np}{2}-1}e^{-\frac{n}{2}%
x}e^{\sum_{j=1}^{r}\theta _{pj}v_{pj}\left( x-S_{p}\right) }\mathrm{d}x,%
\vspace{-0.2in}  \notag
\end{eqnarray}%
where the last equality in (\ref{last eq}) follows from (\ref{LR1}) and
Proposition 2.

The last equality in (\ref{last eq}), (\ref{LR2}) and the fact that$\vspace{%
-0.2in}$%
\begin{equation*}
\int_{p\!-\!\alpha \sqrt{p}}^{p\!+\!\alpha \sqrt{p}}\!\!\!x^{\frac{np}{2}%
-1}e^{-\frac{n}{2}x}e^{\sum_{j=1}^{r}\theta _{pj}v_{pj}\left( x-S_{p}\right)
}\mathrm{d}x\!=\!e^{\sum_{j=1}^{r}-\frac{h_{j}}{2c_{p}}S_{p}}\left( \frac{n}{%
2}\!-\!\sum_{j=1}^{r}\frac{h_{j}}{2c_{p}}\right) ^{-\frac{np}{2}}\!\Gamma
\left( \frac{np}{2}\right) \left( 1\!+\!o(1)\right) \vspace{-0.2in}
\end{equation*}%
imply that$\vspace{-0.2in}$%
\begin{eqnarray*}
L_{p}\!\left( h;\mu _{p}\!\right) &=&\left( 1+o(1)\right) L_{p}\!\left(
h;\lambda _{p}\right) \!e^{\sum_{j=1}^{r}-\frac{h_{j}}{2c_{p}}S_{p}}\left(
1\!-\!\sum_{j=1}^{r}\frac{h_{j}}{nc_{p}}\right) ^{-\frac{np}{2}}\! \\
&=&\left( 1+o(1)\right) L_{p}\!\left( h;\lambda _{p}\right) \!e^{-\frac{%
S_{p}\!-\!p}{2c_{p}}\sum_{j=1}^{r}\!h_{j}+\frac{1}{4c_{p}}\left(
\sum_{j=1}^{r}h_{j}\right) ^{2}},\vspace{-0.2in}
\end{eqnarray*}%
which establishes (\ref{equivalence 2}). The rest of the statements of
Theorem 1 follow from (\ref{equivalence 1}), (\ref{equivalence 2}), and
Lemmas 12 and A2 of OMH.$\square $

\subsection{Proof of Proposition \protect\ref{Proposition3}}

To save space, we only derive the asymptotic power envelope for the
relatively more difficult case of real-valued data and $\mu $-based tests.
According to the Neyman-Pearson lemma, the most powerful test of $h=0$
against the simple alternative $h=\left( h_{1},...,h_{r}\right) $ is the
test which rejects the null when $L_{p}\!\left( h;\mu _{p}\!\right) $ is
larger than a critical value $C.$ It follows from Theorem 1 that, for such a
test to have asymptotic size $\alpha $, $C$ must be$\vspace{-0.2in}$%
\begin{equation}
C=\sqrt{W\left( h\right) }\Phi ^{-1}\left( 1-\alpha \right) +m\left(
h\right) ,\vspace{-0.2in}  \label{ctitical value}
\end{equation}%
where$\vspace{-0.2in}$%
\begin{eqnarray*}
m\left( h\right) &=&\frac{1}{4}\sum_{i,j=1}^{r}\left( \ln \left( 1-\frac{%
h_{i}h_{j}}{c}\right) +\frac{h_{i}h_{j}}{c}\right) \text{ and} \\
W\left( h\right) &=&-\frac{1}{2}\sum_{i,j=1}^{r}\left( \ln \left( 1-\frac{%
h_{i}h_{j}}{c}\right) +\frac{h_{i}h_{j}}{c}\right) .\vspace{-0.2in}
\end{eqnarray*}%
Now, according to Le Cam's third lemma and Theorem 1, under $h=\left(
h_{1},...,h_{r}\right) ,$ $\ln L_{p}\!\left( h;\mu _{p}\!\right) \overset{d}{%
\rightarrow }N\left( m\left( h\right) +W\left( h\right) ,W\left( h\right)
\right) .$ The asymptotic power (\ref{local power mu}) follows.$\square $

\subsection{Invariance issues and Proof of Proposition \protect\ref%
{Proposition4}}

Before turning to the proof of Proposition \ref{Proposition4}, let us
clarify the invariance issues in the problem under study. For basic
definitions (invariant, maximal invariant, etc.), we refer to Chapter 6 of
Lehmann and Romano (2005).

Suppose that $X$ is a $p\times n$ random matrix with $\mathrm{vec}\left(
X\right) \sim N\left( 0,I_{n}\otimes \Sigma \right) $. This model is clearly
invariant under the group $\mathcal{G}_{p},$ acting on $\mathbb{R}^{p\times
n},$ of left-hand multiplications by a $p\times p$ orthogonal matrix $%
x\mapsto Qx,$ $x\in \mathbb{R}^{p\times n},$ $Q\in \mathcal{O}\left(
p\right) ;$ so are the null hypothesis $H_{0}$ and the alternative $H_{1}.$
Letting $m=\min \left( n,p\right) ,$ the $m$-tuple $\lambda \left( X\right)
=\left( \lambda _{1},...,\lambda _{m}\right) $ of non-zero eigenvalues of $%
\frac{1}{n}XX^{\prime }$ is clearly invariant under that group, since $\frac{%
1}{n}xx^{\prime }$ and $\frac{1}{n}\left( Qx\right) \left( Qx\right)
^{\prime }=\frac{1}{n}Qxx^{\prime }Q^{\prime }$ share the same eigenvalues $%
\lambda \left( x\right) $ for any orthogonal matrix $Q$ and any matrix $x\in 
\mathbb{R}^{p\times n}.$ However, $\lambda \left( X\right) $ is not maximal
invariant for $\mathcal{G}_{p},$ as $xx^{\prime }$ and $\left( xP\right)
\left( xP\right) ^{\prime }=xPP^{\prime }x^{\prime }=xx^{\prime },$ where $P$
is an arbitrary $n\times n$ orthogonal matrix, share the same $\lambda
\left( x\right) =\lambda \left( xP\right) $ although, in general, there is
no $p\times p$ orthogonal matrix $Q$ such that $xP=Qx.$

Now, the joint density of the elements of $X$ is%
\begin{equation*}
f_{\Sigma }^{(n)}\left( x\right) =\left( 2\pi \right) ^{-np/2}\left\vert
\Sigma \right\vert ^{-n/2}\exp \left\{ -\frac{1}{2}\mathrm{tr}\left( \Sigma
^{-1}xx^{\prime }\right) \right\} ,\text{ }x\in \mathbb{R}^{p\times n}.
\end{equation*}%
By the factorization theorem, $XX^{\prime }$ is a sufficient statistic, and
it is legitimate to restrict attention to $XX^{\prime }$-measurable
inference procedures. Left-hand orthogonal multiplications $Qx$ of $x$
yields, for $xx^{\prime },$ a transformation of the form $Qxx^{\prime
}Q^{\prime }.$ When $Q$ range over the family $\mathcal{O}_{p}$ of $p\times p
$ orthogonal matrices, those transformations also form a group, $\widetilde{%
\mathcal{G}}_{p},$ say, now acting on the space of $p\times p$ symmetric
positive semidefinite real matrices of rank $m$. Clearly, $\lambda \left(
x\right) $ is maximal invariant for $\widetilde{\mathcal{G}}_{p},$ as $%
xx^{\prime }$ and $yy^{\prime }$ share the same eigenvalues if and only if $%
yy^{\prime }=Qxx^{\prime }Q^{\prime }$ for some $p\times p$ orthogonal
matrix $Q.$

Combining the principles of sufficiency and invariance thus leads to
considering $\lambda $-measurable tests only.

A similar reasoning applies in the case of unspecified $\sigma ^{2},$ with a
larger group combining multiplication by an arbitrary non-zero constant with
the $p\times p$ left orthogonal transformations. Sufficiency and invariance
then lead to restricting attention to $\mu $-measurable tests.

\underline{Proof of Proposition \ref{Proposition4}.}

With the same notation as above, write $T=T\left( X\right) =XX^{\prime }$
for the sufficient statistic. Consider an arbitrary invariant (under the
group $\mathcal{G}_{p}$ of left orthogonal transformations of $\mathbb{R}%
^{p\times n}$) test $\phi \left( X\right) $, and define $\psi \left(
t\right) =E\left( \phi \left( X\right) |T=t\right) $. Then $\psi \left(
T\right) $ is a $T$-measurable test with the same size and power function as 
$\phi \left( X\right) .$ It follows from the proof of Theorem 6.5.3 (i) in
Lehmann and Romano (2005) that $\psi \left( T\right) $ is \textit{almost
invariant}. Moreover, since the conditions of Lemma 6.5.1 (same reference)
hold, this test is invariant under the group $\widetilde{\mathcal{G}}_{p}$
(acting on $T$). Since the ordered $m$-tuple $\lambda _{1},...,\lambda _{m}$
of the eigenvalues of $\frac{1}{n}T=\frac{1}{n}XX^{\prime }$ is maximal
invariant for $\widetilde{\mathcal{G}}_{p},$ and since any invariant
statistic is a measurable function of a maximal invariant one, $\psi \left(
T\right) $ must be $\lambda $-measurable. Hence, $\psi \left( T\right) $ is
a $\lambda $-measurable test and has the same power function as $\phi \left(
X\right) ,$ as was to be shown.

The existence of a $\mu $-measurable test with the same power function as
that of a test $\phi \left( X\right) $ invariant under left orthogonal
transformations and multiplication by non-zero constants is established
similarly.$\square $

\end{document}